\documentclass[draft]{article}
\def\today{1.1.13} 
\usepackage{amsmath,amsfonts,amsthm,amssymb,amscd}
\renewcommand{\Re}{\mathop{\rm Re}\nolimits}
\renewcommand{\Im}{\mathop{\rm Im}\nolimits}

\theoremstyle{plain} \newtheorem{theorem}{Theorem}[section]
\newtheorem{lemma}[theorem]{Lemma}
\newtheorem{proposition}[theorem]{Proposition}
\newtheorem{corollary}[theorem]{Corollary} \theoremstyle{definition}
\newtheorem{definition}[theorem]{Definition} \theoremstyle{remark}
\newtheorem{remark}[theorem]{Remark}

\newcommand{\R}{{\mathbb R}} \newcommand{\U}{{\mathcal U}}
\newcommand{\Hc}{{\mathcal H}} 
\newcommand{\Kc}{{\mathcal K}} 
\newcommand{\Kcs}{\tilde{\mathcal H}}

\newcommand{\Z}{{\mathbb Z}}

\newcommand{\e}{{\rm e}}

\newcommand{\E}{{\mathcal E}}
\newcommand{\F}{{\mathcal F}}

\newcommand{\Tr}{{\mathcal T}}
\newcommand{\Tc}{{\mathcal T}}

\newcommand{\Oc}{{\mathcal O}}
\newcommand{\Wc}{{\mathcal W}}
\newcommand{\Ph}{{\mathcal P}}

\newcommand{\M}{{\mathcal M}}
\newcommand{\I}{{\mathcal I}}
\newcommand{\G}{{\mathcal G}}
\newcommand{\Le}{{\mathcal L}}
\newcommand{\resto}{{\mathcal R}}
\def\im{{\rm i}}

\newcommand{\W}{{\mathcal W}}
\newcommand{\Nc}{{\mathcal N}}
\newcommand{\Sc}{{\mathcal S}}

\newcommand{\V}{{\mathcal V}}
\newcommand{\A}{{\mathcal A}}
\newcommand{\Ca}{{\mathcal C}}

\newcommand{\C}{\mathbb{C}}
\newcommand{\T}{\mathbb{T}}

\def\bq{{\bf q}}

\font\strana=cmti10
\def\lie{\hbox{\strana \char'44}}
\def\uno{{\kern+.3em {\rm 1} \kern -.22em {\rm l}}}

\def\dep#1#2{\frac{\partial{#1}}{\partial p_{#2}}}
\def\di{{\rm d}}
\def\pip{\Pi_p}
\def\depp#1{\dep{\pip}{#1}}
\def\dm{D_t^{-1}}

\def\Phzero{\Ph_0}
\def\Pinr{\Pi_{nr}}

\def\norma#1{\left\| #1\right\|}

\def\sleq{\preceq}

\def\gen{{\rm Gen}}
\def\ac{{\mathcal A}\kern-.7pt\ell\kern-.9pt\mathcal{S}}

\numberwithin{equation}{section}

\begin{document}

\title{Asymptotic stability of ground states in some Hamiltonian
  PDEs with symmetry}
\author{Dario Bambusi}

\date{\today}
\maketitle

\begin{abstract}
We consider a ground state (soliton) of a Hamiltonian PDE. We prove
that if the soliton is orbitally stable, then it is also
asymptotically stable. The main assumptions are transversal
nondegeneracy of the manifold of the ground states, linear dispersion
(in the form of Strichartz estimates) and nonlinear Fermi Golden
Rule. We allow the linearization of the equation at the soliton to
have an arbitrary number of eigenvalues. The theory is tailor made for
the application to the translational invariant NLS in space dimension
3. The proof is based on the extension of some tools of the theory of
Hamiltonian systems (reduction theory, Darboux theorem, normal form)
to the case of systems invariant under a symmetry group with unbounded
generators.
\end{abstract}

\section{Introduction}
\label{1} 

In this paper we study the asymptotic stability of the ground state in
some dispersive Hamiltonian PDEs with symmetry. We will prove that, in
a quite general situation, an orbitally stable ground state is also
asymptotically stable. In order to describe the main result of the
paper we concentrate on the specific model given by the translationally
invariant subcritical NLS in space dimension 3, namely
\begin{align}
\label{NLS1}
 \psi_t=\im\Delta\psi+\im\beta'(\left|\psi\right|^2)\psi\ ,\quad
|\beta^{(k)}(u)|\leq C_k \langle u\rangle^{1+p-k}\ ,\quad
\beta'(0)=0\ .
\end{align}
$p<\frac{2}{3},$ $x\in\R^3$. It is well known that, under suitable
assumptions on $\beta$, such an equation has a family of ground states
which can travel at any velocity and which are orbitally stable (see
e.g. \cite{FroJon} for a review). Consider the linearization of the
NLS at the soliton, and let $L_0$ be the linear operator describing
such a linearized system. Due to the symmetries of the system, zero is
always an eigenvalue of $L_0$ with algebraic multiplicity at least
8. In the case where this is the exact multiplicity of zero and $L_0$
has no other eigenvalues, asymptotic stability was proved in
\cite{BusPer92,Cuc01} (see also \cite{appPer}). Here we tackle the
case where $L_0$ has an arbitrary number of eigenvalues, disjoint from
the essential spectrum, and prove that, assuming a suitable version of
the Fermi Golden Rule (FGR), the ground state is (orbitally)
asymptotically stable. We recall that the importance of the FGR in
nonlinear PDEs was understood by Sigal \cite{Sig93} and shown to have
a crucial role in the study of asymptotic stability in
\cite{SofW99}. Similar conditions have been used and generalized by
many authors. The FGR that we use here is a generalization of that of
\cite{GW08} (see also \cite{BC09,Cuc09}).
 
\vskip10pt

 The present paper is a direct development of \cite{BC09} and
 \cite{Cuc09}, which in turn are strongly related to
 \cite{GS07,GSNT04,CucMiz,GW08}. We recall that in \cite{BC09}
 Hamiltonian and dispersive techniques were used to prove that the
 empty state of the nonlinear Klein Gordon equation is asymptotically
 stable even in the presence of discrete spectrum of the linearized
 system. Then \cite{Cuc09} extended the techniques 
 of \cite{BC09} to the study of the asymptotic stability of the ground
 state in the NLS with a potential.

The main novelty of the present paper is that we deal here with the
translational invariant case. The new difficulty one has to tackle is
related to the fact that the group of the translations
$\psi(.)\mapsto\psi(.-t{\bf e}_i)$ is generated by $-\partial_{x_i}$
which is an unbounded operator: it turns out that this obliges to use
non smooth maps in order to do some steps of the proof. To overcome
this problem we introduce and study a suitable class of maps, that we call
``almost smoothing perturbation of the identity'' (see in particular sect.\ref{s.darboux}). We use
them to develop Hamiltonian reduction theory, Darboux theorem and also
canonical perturbation theory.

The fact that the generator of the translations is not smooth causes some
difficulties also in the use of Strichartz estimates, but such
difficulties were already overcome by Beceanu \cite{Bec11} and by Perelman \cite{appPer}, so we simply apply their method to our case.

\vskip10pt We now describe the proof.  First, we use Marsden Weinstein
reduction procedure in order to deal with the symmetries. In order to
overcome the problems related to the fact that the generators of the
symmetry group are unbounded, we fix a concrete local model for the
reduced manifold and work in it. The local model is a submanifold
contained in the level surface of the integrals of motion. The
restriction of the Hamiltonian and of the symplectic form to such a
submanifold give rise to the Hamiltonian system one has to study. The
advantage of such an approach is that the ground state appears as a
minimum of the Hamiltonian, so one is reduced to study the asymptotic
stability of an elliptic equilibrium, a problem close to that studied 
in \cite{BC09}. However the application of the methods of
\cite{BC09,Cuc09} to the present case is far from trivial, since the
restriction of the symplectic form to the submanifold turns out to be
in noncanonical form, and to have non smooth coefficients (some
``derivatives'' appear).  So, we proceed by first proving a suitable
version of the Darboux theorem which reduces the symplectic form to
the canonical one.  This requires the use of non smooth
transformations. We point out that a key ingredient of our
developments is that the ground state is a Schwartz function, and this
allows to proceed by systematically moving derivatives from the
unknown function to the ground state.

Then we study the structure of the Hamiltonian in the Darboux
coordinates and prove that it has a precise (and quite simple)
form. Subsequently, following \cite{BC09,Cuc09}, we develop a suitable
version of normal form theory in order to extract the essential part
of the coupling between the discrete modes and the continuous
ones. Here we greatly simplify the theory of \cite{BC09,Cuc09}. In
particular we think that we succeeded in developing such a theory
under minimal assumptions. We also point out that in the present case
the canonical transformations putting the system in normal form are
not smooth, but again almost smoothing perturbation of the identity.

Finally, following the scheme of \cite{GSNT04,CucMiz,BC09,Cuc09}, we
use Strichartz estimates in order to prove that there is dispersion,
and that the energy in the discrete degrees of freedom goes to zero as
$t\to\infty$. As we already remarked there are some difficulties in
the linear theory, difficulties that we overcome using the methods of
\cite{appPer}. In this part, we made an effort to point out the
properties that the nonlinearity has to fulfill in order to ensure the
result. Thus we hope to have proved a result which can be simply
adapted to different models.

\vskip5pt

We now discuss more in detail the relation with the paper
\cite{Cuc09}. In \cite{Cuc09} Cuccagna studied the case of NLS with a
potential and proved a result similar to the present one. Here we
generalize Cuccagna's result in several aspects. The first one is that
we allow the system to have symmetry groups with more than one
dimension, but the main improvement we get consists of the fact that we
allow the symmetries to be generated by unbounded operators (as
discussed above). Furthermore we work in an abstract framework.

Finally, we work here on the
reduced system (according to Marsden-Wein\-stein theory), but we think
that all the arguments developed in such a context could be reproduced
also working in the original phase space. We also expect that the same
(maybe more) difficulties will appear also when working in the original
phase space.  
\vskip5pt 

\vskip 5pt Three days before the first version of this paper was
posted in Arxiv, the paper \cite{CucMov} was also posted there. The
paper \cite{CucMov} deals exactly with the same problem. The result of
\cite{CucMov} is very close to the present one, but weaker:
the result of such a paper is valid only for initial data of Schwartz
class, while the control of the difference between the soliton and the
solution is obtained in energy norm, and no decay rate is
provided. Such a kind of conclusions is usual for initial data in the
energy space, while the typical result valid for solutions
corresponding to initial data decaying in space also controls the rate
of decay of the solution to the ground state.  On the contrary, in the
present paper we give a result valid for any initial datum of finite energy
(and of course we do not deduce a decay rate).

A further difference between the two papers is that, here
a large part of the proof is developed in an abstract framework, thus
we expect our result to be simply applicable also to different
systems. We are not aware of other papers in the domain of asymptotic
stability in dispersive Hamiltonian PDEs in which the proof is
developed in an abstract framework.

Our proof is also much simpler than that of \cite{CucMov}, indeed in
order to generate the flow of the transformation introducing Darboux
coordinates (and the transformations putting the system in normal
form) we use a technique coming from the theory of semilinear PDEs,
while \cite{CucMov} uses techniques coming from quasilinear PDEs. 

A further difference, is that we work using Marsden Weinstein
reduction, while \cite{CucMov} works in the original phase space.

\vskip5pt

The paper is organized as follows: in sect. \ref{StateNLS} we state
our main result for the NLS; in sect. \ref{2} we set up the abstract
framework in which we work and state and prove the Darboux theorem
mentioned above; in sect. \ref{hami} we study the form of the
Hamiltonian in the Darboux coordinates. 
In sect. \ref{Normalform} we use canonical perturbation theory in
order to decouple as far as possible the discrete degrees of freedom
from the continuous ones; in sect. \ref{DisNLS} we prove that the
variables corresponding to the continuous spectrum decay dispersively
and the variables corresponding to the discrete spectrum decay at
zero; here the main abstract theorem \ref{main.dis} is stated and
proved; in sect. \ref{NLS} we apply the abstract theory to the
NLS. In the first Appendix we prove that the dynamics of the reduced
system, while in the second one, we reproduce Perelman's Lemma
on the dispersion of the linear system. 

\noindent
{\it Acknowledgments}. First, I would like to warmly thank Galina
Perelman for sending me her notes on asymptotic stability of solitons
in energy space. During the preparation of this paper I benefit of the
constructive criticism and of the suggestions by many persons. In
particular I would like to thank N. Burq, P. Gerard, S. Gustafson,
T. Kappeler, E. Terraneo. In the second version of \cite{CucMov} some
criticisms are raised on a previous version of the present
paper. The analysis stimulated by such criticisms led me to a
considerable simplification of the proof.

\section{Asymptotic stability in NLS}\label{StateNLS}

We state here our result on the NLS eq. \eqref{NLS1}. We assume
\begin{itemize}
\item[(H1)] There exists an open interval $\I\subset \R$ such that,
  for ${\E}\in\I$ the
  equation
\begin{equation}
\label{ground.NLS.0}
-\Delta b_{\E}-\beta'(b_{\E}^2)b_{\E}+{\E}b_{\E}=0\ ,
\end{equation}
admits a $C^{\infty}$ family of positive, radially symmetric functions
$b_{\E}$ belonging to the Schwartz space.
\item[(H2)] One has $\displaystyle{
\frac{ d}{d{\E}}\norma{b_{\E}}^2_{L^2}>0,\ {\E}\in\I}$.
\end{itemize}
\noindent
Then one can construct traveling solitons, which are 
solutions of \eqref{NLS1} of the form 
\begin{equation}
\label{sol1}
\psi(x,t)=e^{-\im\left(\E-\frac{|v|^2}{4}\right)t} \e^{-\im\frac{v\cdot
    x}{2}}b_{\E}(x-vt)\ .
\end{equation} 

\begin{itemize}
\item[(H3)] Consider the operators
\begin{equation}
\label{sta.2}
A_+:=-\Delta+{\E}-\beta'(b_{\E}^2)\ ,\quad A_-:=-\Delta+\E
-\beta'(b_\E ^2)-2\beta''(b_\E ^2)b^2_\E\ ,
\end{equation}
then the Kernel of the operator $A_+$ is generated by $b_{\E}$ and
the Kernel of the operator $A_-$ is generated $\partial _jb_\E$,
$j=1,2,3$. 
\end{itemize}

\begin{remark}
\label{stab}
Under the above assumptions the solutions \eqref{sol1} are orbitally
stable (see e.g. \cite{FroJon}).
\end{remark}
In order to state the assumptions on the linearization at the soliton
insert the following Ansatz in the equations
\begin{equation}
\label{stab.3}
\psi(x,t)=e^{-\im\left({\E}-\frac{|v|^2}{4}\right)t} e^{-\im\frac{v\cdot
    x}{2}}(b_{\E}(x-vt)+\chi(x-vt))\ ,
\end{equation} 
and linearize the so obtained equation in $\chi$. Then one gets an
equation of the form $\dot\chi= L_0\chi$ with a suitable $L_0$.  It
can be easily proved that the essential spectrum of $L_0$ is
$\bigcup_{\pm}\pm \im [{\E},+\infty)$ and that 0 is always an
  eigenvalue. The rest of the spectrum consists of purely imaginary
  eigenvalues $\pm\im \omega_j$, that we order as follows
  $0<\omega_1\leq\omega_2\leq...\leq \omega_K$. We assume that

\begin{itemize}
\item[(H4)] $\omega_K<{\E}$. Furthermore, let $r_t$ be the smallest
  integer number such that $r_t\omega_1>{\E}$, then we assume
$\omega\cdot k\not={\E}$, $ \forall k\in\Z^K\ :\ |k|\leq 2r_t.$

\item[(H5)] $\pm \im {\E}$ are not resonances of $L_0$.  
\item[(H6)] The Fermi Golden Rule \eqref{m.34} holds.
\end{itemize}

The main theorem we are now going to state refers to initial data
$\psi_0$ which are sufficiently close to a ground state. In its
statement we denote by $\epsilon$ the quantity below
\begin{equation}
\label{eps.1}
\epsilon:=\inf_{q^4_0\in\R,\bq_0\in\R^3, v_0\in\R^3,{\E}_0\in\I}\norma{\psi_0-
  \e^{-\im q^4_0}\e^{-i\frac{v_0\cdot x}{2}}b_{{\E}_0}(x-\bq_0)}_{H^1}
\end{equation}

\begin{theorem}
\label{main.NLS}
Assume $\epsilon$ is small enough, then there exist $C^1$ functions
$${\E}(t),v(t),q^4(t), {\bf q}(t),y^4(t),{\bf y}(t)\ , $$ and
$\psi_+\in H^1$ such that the solution $\psi(t)$ with initial datum
$\psi_0$ admits the decomposition
\begin{equation}
\label{main.nls.1}
\psi(x,t)=\e^{-\im q^4(t)}\e^{-\im\frac{v(t)\cdot
    x}{2}}b_{{\E}(t)}(x-\bq(t))+\e^{-\im y^4(t)} \chi(x-{\bf y}(t),t)
\end{equation}
and
\begin{equation}
\label{main.nls.2}
\lim_{t\to+\infty}\norma{\chi(t)-\e^{\im t\Delta}\psi_+}_{H^1}=0\ .
\end{equation}
Furthermore the functions ${\E}(t),v(t), \dot q^4(t),\dot y^4(t),\dot
\bq(t), \dot {\bf y}(t)$ admit a limit as $t\to+\infty$. 
\end{theorem}
The rest of the paper is devoted to the proof of an abstract version
of this theorem.

\section{General framework and the Darboux theorem}\label{2}

Consider a scale of Hilbert spaces $\Hc\equiv\{\Hc^k\}_{k\in\Z}$,
$\Hc^{k+1}\hookrightarrow \Hc^k$ continuously. The scalar product in
$\Hc^0$ will be denoted by $\langle.;.\rangle $; such a scalar product
is also the pairing between $\Hc^k$ and $\Hc^{-k}$. We will denote
$\Hc^\infty:=\cap_{k}\Hc^k$, and $\Hc^{-\infty}:=\cup_{k}\Hc^k$.  Let
${E}:\Hc^k\to \Hc^{k}$, $\forall k$ be a linear continuous operator
skew-symmetric with respect to $\langle.;.\rangle$. Assume it is
continuously invertible. Let $J:\Hc^k\to \Hc^{k}$ be its inverse
(Poisson tensor). We endow the scale by the symplectic form
$\omega(U_1,U_2):=\langle EU_1;U_2\rangle$, then the Hamiltonian
vector field $X_H$ of a function $H$ is defined by $X_H=J\nabla
H$,where $\nabla H$ is the gradient with respect to the scalar product
of $\Hc^0$.

\begin{remark}
\label{2indici}
In the application to dispersive equations one has to deal with
weighted Sobolev space $H^{k_1,k_2}$, which are labeled by a couple of
indexes. All what follows holds also in such a situation provided one
defines the notation $(k_1,k_2)>(l_1,l_2)$ by  $k_1>l_1$ and $k_2\geq l_2$.
\end{remark}

For $j=1,...,n$, let $A\relax_j: \Hc^k\to \Hc^{k-d_j}$, $\forall
k\in\Z$ and some $d_j\geq 0$, be $n$ bounded selfadjoint (with respect
to $\langle.;.\rangle$) linear operators, and consider the
Hamiltonian functions $\Ph\relax_j(u):=\langle A\relax_j u;u\rangle/2
$. Then $X_{\Ph\relax_j}=JA\relax_j$ generates a flow in $\Hc^0$
denoted by $\e^{t JA\relax_j}$.
\begin{remark}
\label{2.ind}
In the case of multiple indexes the index $d_j$ represents the loss of
smoothness and always acts only on the first index, namely one has
$A_jH^{k_1,k_2}\subset H^{k_1-d_j,k_2}$.
\end{remark}
\begin{remark}
\label{rem.sym.prima}
The operators $JA_j$ will play the role of the generators of the
symmetries of the Hamiltonian system we will study. Correspondingly
the functions $\Ph_j$ will be integrals of motion.
\end{remark}

We denote $d_A:=\max_{j=1,...,n}d_j$.  For $i,j=1,...,n$ we assume
that, on $\Hc^{\infty}$ one has

\begin{itemize}

\item[(S1)] $[A_j,E]=0$,

\item[(S2)] $A\relax_iJA\relax_j=A\relax_jJA\relax_i$ which implies
  $\left\{\Ph\relax_j,\Ph\relax_i\right\}=0=\langle
  A\relax_ju;JA\relax_iu\rangle$.
\item[(S3)] For any $t\in\R$ the map $\e^{tJA_j}$ leaves invariant
  $\Hc^k$ for all $k$ large enough.

\end{itemize}

Let $A_0$ be a linear operator with the same properties of the
$A_j$'s. Assume $d_0\geq d_A$.  The Hamiltonian we will study has the
form
\begin{equation}
\label{e.ham}
H(u)=\Phzero\null(u)+H_P(u)\ ,\quad \Phzero\null(u)=\frac{1}{2}\langle
u; A_0 u\rangle\ ,
\end{equation}
where $H_P$ is a nonlinear term on which we assume
\begin{itemize}
\item[(P1)] There exists $k_0$ and an open neighborhood of zero
  $\U^{k_0}\subset \Hc^{k_0}$ such that $H_P\in
  C^{\infty}(\U^{k_0},\R)$.
\end{itemize}

We also assume that (on $\Hc^{\infty}$)

\begin{itemize}
\item[(S4)] $H_P$ and $\Ph_0$ Poisson commutes with each one of the functions
$\Ph\relax_j$:
\begin{equation}
\label{e.ham.2}
\left\{\Ph_0;\Ph\relax_j\right\}=\left\{
H_P;\Ph\relax_j\right\}=0\ ,\quad j=1,...,n
\end{equation}
\end{itemize}

We are interested in bound states $\eta$, namely in phase points such
that $u(t):=\e^{t \lambda^j JA\relax_j}\eta$ is a solution of the
Hamilton equations of $H$. {\it Here and below we use Einstein
  notation according to which sum over repeated indexes is
  understood.} The indexes will always run between $1$ and $n$.  Then
$\eta$ has to fulfill the equation
\begin{equation}
\label{e.ham.4}
A_0\eta+\nabla H_P(\eta)-\lambda^j A\relax_j\eta=0\ .
\end{equation} 

We assume
\begin{itemize}
\item[(B1)] There exists an open set $I\subset \R^n$ and a
  $C^{\infty}$ map
$$I\ni p \mapsto (\eta_p,\lambda(p))\in
  \Hc^\infty\times \R^n\ ,$$ s.t.
  $(\eta_p,\lambda(p))$ fulfills equation \eqref{e.ham.4}. Furthermore
  the map $p\to \lambda$ is 1 to 1.

\item[(B2)] For any fixed $p\in I$, the set
  $\Ca:=\bigcup_{q\in\R^n}\e^{q^jJA_j}\eta_p$ is a smooth $n$
  dimensional submanifold of $\Hc^\infty$.
\item[(B3)] The manifold $\bigcup_{p\in I}\eta_p$ is isotropic, namely
  the symplectic form $\omega$ vanishes on its tangent space.
\end{itemize}

By (B1) it is possible to normalize the values of $p_j$ in such a way
that $\Ph_j(\eta_p)=p_j$, {\sl Form now on we will always assume such
  a condition to be satisfied.} In particular it follows (by deriving
with respect to $p_k$), that 
\begin{equation}
\label{nor.5.2}
\left\langle A_j \dep{\eta_p}{k};\eta_p \right\rangle=\delta_j^k\ .
\end{equation}

\begin{remark}
\label{r.31}
By the proof of Arnold Liouville's theorem, the manifold $\Ca$ of
hypothesis (B2) is diffeomorphic to $\T^k\times \R^{n-k}$, where
$\T=\R/2\pi\Z$.
\end{remark}

Consider the symplectic manifold $$\Tr:=\bigcup_{q\in \R^n,p\in
  I}\e^{q^jJA_j}\eta_p\ ,$$ namely the manifold of bound states; its
tangent space is given by
\begin{equation}
\label{i.4}
T_{\eta_p}\Tr:={\rm span}\left\{ JA\relax_j\eta_p,\frac{\partial
  \eta_p}{\partial p_j}\right\}
\end{equation} 
and its symplectic orthogonal $T^\omega _{\eta_p}\Tr$ is given by
\begin{align}
\nonumber
& T^\omega _{\eta_p}\Tr 
\\
\nonumber
 &=\left\{ U\in
\Hc^{-\infty}\ :\omega\left(JA\relax_j\eta_p;U\right)=\left\langle
A\relax_j\eta_p;U\right\rangle=
   \omega\left(\dep{\eta_p}j;U\right)=\langle E\dep{\eta_p}j;U\rangle=0
\right\} 
\end{align}

\begin{lemma}
\label{deco}
One has $\Hc^{-\infty}=T^\omega_{\eta_p}\Tr\oplus T_{\eta_p}\Tr$. 
Explicitly the decomposition of a vector $U\in \Hc^{-\infty}$ is given by
\begin{equation}
\label{i.5a}
U=P\relax_j\dep {\eta_p}j+Q^jJA\relax_j\eta_p+\Phi_p\ ,
\end{equation}
with 
\begin{equation}
\label{i.5b}
Q^j=-\langle E\dep{\eta_p}j;U\rangle\ ,\quad \quad P\relax_j=\langle
A\relax_j\eta_p;U\rangle 
\end{equation}
and $\Phi_p\in T^{\omega}_{\eta_p}\Tr$ given by 
\begin{equation}
\label{i.5cr}
\Phi_p=\Pi_pU:=U-\langle A\relax_j\eta_p;U\rangle\dep{\eta_p}j+ \langle E\dep 
    {\eta_p}j;U \rangle JA\relax_j\eta_p\ .
\end{equation}
\end{lemma}
\proof The first of \eqref{i.5b} is obtained taking the scalar product
of \eqref{i.5a} with $-E\dep{\eta_p}j$, exploiting \eqref{nor.5.2}
and  
\begin{equation}
\label{iso}
\langle E\dep{\eta_p}j,\dep{\eta_p}k\rangle= 0
\end{equation}
which is equivalent to (B3).  Taking the scalar product of
\eqref{i.5a} with $A\relax_j\eta_p$ we get the second of
\eqref{i.5b}. Then \eqref{i.5cr} immediately follows.\qed

\begin{remark}
\label{smo.pi}
A key point in all the developments of the paper is that the projector
$\Pi_p$ defined by \eqref{i.5cr} is a smoothing perturbation of the
identity, namely $\uno-\Pi_p\in C^{\infty}(I,B(\Hc^{-k},\Hc^{l}))$,
$\forall k,l $, where $B(\Hc^{-k},\Hc^{l}) $ is the space of
bounded operators from $\Hc^{-k}$ to $\Hc^{l}$. In
particular one has that
$\depp {i}\in C^{\infty}(I,B(\Hc^{-k},\Hc^{l}))$.
\end{remark}
 An explicit computation shows that the adjoint of $\pip$ is given by 
\begin{equation}
\label{i.5f}
\Pi_p^*U:=U-\langle\dep{\eta_p}j;U\rangle A\relax_j\eta_p+\langle
JA\relax_j\eta_p,U\rangle E\dep{\eta_p}j\ .
\end{equation}

Some useful formulae are collected below
\begin{align}
\label{r.2.1}
E\Pi_p=\Pi^*_pE\ ,\quad J\Pi_p^*=\Pi_pJ\ ,\quad
\dep{\Pi_p}j=\dep{\Pi_p^2}j=\Pi_p\dep{\Pi_p}j+\dep{\Pi_p}j\Pi_p \ .
\\
\label{r.2.3}
\Pi_p\dep{\Pi_p}j\Pi_p=0\ ,\quad
\left(\dep{\Pi_p}j\right)^*=\dep{\Pi_p^*}j\ ,\quad
E\dep{\Pi_p}j=\dep{\Pi_p^*}j E\ .
\end{align}

{\bf In the following we will work locally close to a particular value
$p_0\in I$. Thus we fix} it and define
\begin{equation}
\label{i.7aa}
\V^k:= \Pi_{p_0}\Hc^{k}
\end{equation}
which we endow by the topology of $\Hc^k$. 
\begin{remark}\label{r.2}
For any positive $k,l$, one has
\begin{align}
\label{i.22}
\norma{(\Pi_p\Pi_{p'}-\Pi_{p'})u}_{\Hc^k}\leq
C_{k,l}\left|p-p'\right|\norma u_{\Hc^{-l}}\ ,  
\end{align}
and, by the first of \eqref{r.2.3}, one has
\begin{equation}
\label{i.23a}
\norma{\Pi_p\dep{\Pi_p}j\phi}_{\Hc^k}\leq
C_{kl}\left|p-p_0\right|\norma{\phi} _{\Hc^{-l}}\ .
\end{equation}
\end{remark}

\begin{remark}
\label{r.5.1}
Consider the operator $\pip:\V^{-\infty}\to \pip \Hc^{-\infty}$; it
has the structure $\uno+(\pip-\Pi_{p_0})$, and one has
\begin{equation}
\label{r.e.5.1}
\norma{(\pip-\Pi_{p_0})\phi}_{\Hc^k}\leq
C\left|p-p_0\right|\norma\phi_{\Hc^{-l}}\ .
\end{equation} 
Thus, by Neumann formula the inverse $\widetilde {\pip}^{-1}$ of
$\pip$ has the form $\widetilde \pip^{-1}=\uno+S$ with $S$ fulfilling
\eqref{r.e.5.1}.  
\end{remark}

\subsection{Reduced manifold}\label{red}

We introduce now the reduced symplectic manifold obtained by
exploiting the symmetry. In the standard case where
the generators of the symmetry group are smooth (i.e. $d_j=0$) the
construction is standard and goes as follows.

Fix $p_0\in I$ as above and define a surface
$\Sc=\{u\ :\ \Ph\relax_j(u)=p_{0j}\}$, then pass to the quotient with
respect to the group action of $\R^n$ on $\Sc$ defined by
$(q,u)\mapsto \e^{q^jJA\relax_j}u $, obtaining the reduced phase space
$\M$. A local model of $\M$ close to $\eta_{p_0}$ is obtained by
taking a codimension $n$ submanifold of $\Sc$ transversal to the orbit
of the group.  Here we proceed the other way round: we choose a
submanifold $\M\subset\Sc$ of codimension $n$, transversal to the
orbit of the group at $\eta_{p_0}$, and we study the Hamiltonian
system obtained by restricting the Hamiltonian to $\M$.

Fix some $k\geq 0$, and consider the map
\begin{equation}
\label{i.7a}
I\times\V^k\ni (p,\phi)\mapsto
i_{0}(p,\phi):=\eta_p+\Pi_p\phi \ ;
\end{equation}
we will use the implicit function theorem (see lemma \ref{impl.1})
in order to compute $p\relax_j=p\relax_j(\phi)$ in such a way that the
image of the map
\begin{equation}
\label{i.7b}
\V^k\ni\phi\mapsto i(\phi):=\eta_{p(\phi)}+\Pi_{p(\phi)}\phi\subset
\Sc \ ,
\end{equation}
is the wanted local model of $\M$, and $i$ is a local coordinate
system in it. In studying this map we will use a class of maps which
will play a fundamental role in the whole paper. First we introduce a
suitable notion of smooth and regularizing map between scales of
Hilbert spaces. 

Let $\Hc\equiv\left\{\Hc^k\right\}$ and
$\Kcs\equiv\left\{\Kcs^l\right\}$ be two scales of Hilbert spaces, then
we give the following definitions

\begin{definition}
\label{def.smo.new}
A map $f$ will be said to be of class $\ac(\Hc,\Kcs)$ if,
$\forall r,l\geq 0$ there exist $k$ and an open neighborhood of the origin
$\U_{rlk}\subset\Hc^k$ such that
\begin{equation}
\label{de.eq.1}
f\in C^r(\U_{rlk};\Kcs^l)\ .
\end{equation}
Sometimes such maps will be called {\it almost smooth}.
\end{definition}
We will use the same notation also when one of the two scales, or both
are composed by a single space.

\begin{remark}
\label{HPAS}
If $\Tr\in\ac(\Hc,\Hc)$ then, by (P1), $H_P\circ \Tr\in\ac(\Hc,\R)$.
\end{remark}

\begin{definition}
\label{def.smo.new.1}
A map $f$ will be said to be regularizing or of class $C_R(\Hc,\Kcs)$ if,
$\forall r,l,k\geq0$ there exists an open neighborhood of the origin
$\U_{rlk}\subset\Hc^{-k}$ such that
\begin{equation}
\label{de.eq.11}
f\in C^r(\U_{rlk};\Kcs^l)\ .
\end{equation}
\end{definition}

In the following the width of open sets does not play any role so
we will avoid to specify it. In particular we will often consider maps from
an Hilbert space (typically $\R^n\times \V^{k}$) to some other
space, {\it by this we {\bf always} mean a map defined in an
  open neighborhood of the origin.}

\begin{definition}
\label{d.1}
For $i,j\geq 0$, a map $S$ will be said to be of class $\Sc_j^i$ if
there exists a regularizing map $\tilde S\in C_R(\Kcs,\Hc)$, $\Kcs:=
\R^n\oplus \V$ such that $S(\phi)=\tilde S(\Ph(\phi),\phi)$ for
sufficiently small $\phi\in\V^{d_A/2}$, and the map $\tilde S$
fulfills, $\forall m,k\geq 0$,
\begin{equation}
\label{cla.1.1}
\norma{\tilde S(N,\phi)}_{\Hc^m}\leq C_{mk}
|N|^i\norma{\phi}^j_{\Hc^{-k}}\ ,
\end{equation} 
$\forall (N,\phi)$ in some neighborhood of the origin in
$\R^n\times\V^{-k}$ (which can depends on $m$).
\end{definition}

In the case of maps taking values in $\R^n$ we give an analogous
definition.
\begin{definition}
\label{d.1.11}
For $i,j\geq 0$, a map $R$ will be said to be of class $\resto_j^i$ if
there exists a regularizing map $\tilde R\in C_R(\Kcs,\R^n)$, $\Kcs:=
\R^n\oplus \V$, such that $R(\phi)=\tilde R(\Ph(\phi),\phi)$ for
sufficiently small $\phi\in\V^{d_A/2}$, and $\forall k\geq 0$ the map
$\tilde R$ fulfills
\begin{equation}
\label{cla.1.11}
\norma{\tilde R(N,\phi)}\leq C_{k}
|N|^i\norma{\phi}^j_{\Hc^{-k}}\ , 
\end{equation} 
$\forall (N,\phi)$ in some neighborhood of the origin in
$\R^n\times\V^{-k}$.
\end{definition}

{\sl The functions belonging to the classes of definitions \ref{d.1}
  and \ref{d.1.11} will be called {\it smoothing}. }

In the following we will identify a smoothing function $S$ (or $R$)
with the corresponding function $\tilde S$ (or $\tilde R$).  Most of
the times functions of class $\Sc^k_l$ ($\resto_k^l$ resp.) will be
denoted by $S^k_l$ ($R^k_l$ resp.). Furthermore, since the only
relevant property of such functions are given by the inequalities
\eqref{cla.1.1} and \eqref{cla.1.11} we will use the same notation for
different smoothing functions. For example we will meet equalities of
the form
\begin{equation}
\label{exem}
S^1_1+S^1_2=S^1_1
\end{equation}
where obviously the function $S^1_1$ at r.h.s. is different from that
at l.h.s.

Finally, we always consider functions and vector fields as functions
of $N,\phi$, with the idea that, at the end of the procedure we will
put $N_j=\Ph_j(\phi)$.

\begin{lemma}
\label{impl.1}
There exists a smoothing map $p\in\resto^0_0 $ with the following
properties
\begin{itemize}
\item[(1)] For any $j=1,...,n$, and for $\phi\in \V^{d_A/2}$, one has 
$$\Ph\relax_j(\eta_{p(\Ph(\phi),\phi)}+\Pi_{p(\Ph(\phi),\phi)}\phi)=p_{0j}\ ; $$
\item[(2)] there exist $R_2^1\in\resto^1_2$ s.t. $
p=p_0-N+R_2^1(N,\phi)$;
\item[(3)] Define the matrix $M=M(N,\phi)$, by $
  \displaystyle{(M^{-1})_{jk}=\delta_{jk}+\langle\Pi_p\phi;A\relax_j\dep{\Pi_p}k\phi\rangle}
  $ (evaluated at $p=p(N,\phi)$), then the gradient of
  $p(\Ph(\phi),\phi)$ is given by
\begin{align}
\label{smo.01}
\nabla p\relax_j= -\sum_kM_{jk}\Pi^*_{p_0} A\relax_k\phi\ .
\end{align}
\end{itemize} 
\end{lemma}
\proof First remark that one has 
\begin{align}
\label{smo.2}
\Ph\relax_{j}(\eta_p+\Pi_p\phi)=p\relax_j+\Ph\relax_j(\Pi_p\phi)
=p\relax_j+\Ph\relax_j(\phi)+ V_j(p,\phi)\ ,
\end{align}
where $V_j(p,\phi):=\Ph_j(\Pi_p\phi)-\Ph_j(\phi)$ extends to a
regularizing map $(p,\phi)\mapsto V_j(p,\phi)$,
$V\in C_R(\R^n\oplus\V,\R^n)$ which $\forall k$ fulfills
\begin{equation}
\label{smo.3}
|V_j(p,\phi)|\leq C_k
\left|p-p_0\right|\norma{\phi}_{\Hc^{-k}}^2\ .\quad 
\end{equation}
We fix $k$ and apply the implicit function theorem to the system
of equations
\begin{equation}
\label{smo.31}
0=F_j(p,N,\phi):=p\relax_j+N\relax_j+V_j(p,\phi)-p_{0j}\ 
\end{equation}
in $(p,N,\phi)\in\R^n\times\R^n\times \V^{-k}$. Using the definition
of $V_j$, one gets
$$ \frac{\partial F_j}{\partial p_k}=\delta^k_j+
\langle\Pi_p\phi;A_j\depp k\phi\rangle= (M^{-1})_{jk}\ .
$$ Since $\dep{F\relax_j}k\equiv (M^{-1})_{jk}$ is invertible, the
implicit function theorem ensures the existence of a smooth function
$p=p(N,\phi)$ from $\R^n\times \Hc^{-k}$ to $\R^n$ solving
\eqref{smo.31}.  Then the estimate ensuring $p-(p_0-N)\in\resto_2^1$
follows from the fact that $p(p_0,\phi)-(p_0-N)=0$ and from the
computation of the differential of $p$ with respect to $\phi$, which
gives
$
\di_\phi p(N,\phi)= -M\di_\phi V
$, which in turn shows that $\di_\phi p(N,0)=0$ since $V$ is quadratic
in $\phi$. Equation \eqref{smo.01} is an immediate consequence of the
formula for the derivative of the implicit function applied directly
to $\eqref{smo.2}=p_0$.  \qed

We are now going to study the correspondence between the dynamics in
the reduced phase space $\V^{k_0}$, (for some $k_0$) and the dynamics
of the complete system. We endow $\V^0$ by the symplectic form
$\Omega:=i^*\omega$ (pull back) and consider an invariant Hamiltonian
function $H:\Hc^{k_0}\to\R$, namely a function with the property that
$H(e^{q^jJA\relax_j}u)=H(u)$.  We will denote by $H_r:=i^*H$ the
corresponding reduced Hamiltonian.

\begin{remark}
\label{hr.1}
By the smoothness of the map $p\mapsto \eta_p$, there exists a map
$S^1_1\in\Sc^1_1$ such that
$H_r(\phi)=H(\eta_{p(N,\phi)}+\Pi_{p(N,\phi)}\phi)=H(\eta_{p_0-N}+\phi+S^1_1)$,
a formula which will be useful in the following.
\end{remark}

The function $H_r$ defines a Hamiltonian system on $\V^{k_0}$. We will
denote by $X_{H_r}$ the corresponding Hamiltonian vector field. Denote
also by $X_H$ the Hamiltonian vector field of $H$ in the original
phase space. Before stating the theorem on the correspondence of the
solutions we specify what we mean by solution.

\begin{definition}
\label{sol}
Let $k$ be fixed. A function $u\in C^0([0,T];\Hc^k)$ will be said to be
a solution of $\dot u=X(u)$, if there exists a sequence of functions
$u_l\in C^1([0,T];\Hc^k)$ which fulfill the equation and converge to $u$
in $C^0([0,T];\Hc^k)$. 
\end{definition}

\begin{theorem}
\label{t.red.1}
Assume that $X_{H_r}$ defines a local flow in $\V^l$ for some $l\geq
0$. Assume that such a flow leaves invariant $\V^k$ for some $k$ large
enough. Let $u_0:=e^{q^j_0JA_j}i(\phi_0)$, be an initial datum with
$\phi_0\in\V^l$. Consider the solution $\phi(t)\in\V^l$ of the Cauchy
problem
$\dot \phi=X_{H_r}(\phi)$, $\phi(0)=\phi_0$.
Then there exist $C^1$ functions $q^j(t)$ such that
\begin{equation}
\label{t.red.e}
u(t):=e^{q^j(t)JA\relax_j}i(\phi(t))
\end{equation}
is a solution of 
$\dot u=X_{H}(u)$
with initial datum $u_0$.
Viceversa, if $X_H$ generates a local flow for initial data close to
$\Ca$ and such a flow leaves invariant $\Hc^k$ for some $k$ large enough,
then any solution of the original system admits the representation
\eqref{t.red.e}, with $\phi(t)$ a solution of the reduced system.
\end{theorem}

The proof is obtained more or less as in the standard way (see
e.g. \cite{Schmid}), however one has to verify that all of the
formulae that are used keep a meaning also in the present non smooth
case, and this is quite delicate. For this reason the proof is
deferred to the appendix \ref{red.proof}.

\subsection{Almost smoothing perturbations of the identity and the
  Darboux theorem}\label{s.darboux} 

Denote $\Omega:=i^*\omega$. By construction it is clear that
\begin{equation}
\label{i.8}
\Omega\big\vert_{\phi=0}\left(\Phi_1;\Phi_2
\right)= \langle E\Phi_1;\Phi_2 \rangle\ .
\end{equation}
We will transform the coordinates in order to obtain that in a whole
neighborhood of 0 the symplectic form takes the form \eqref{i.8}.
The coordinate changes we will use are not smooth (this would be
impossible, since, due to our construction, the symplectic form
$\Omega$ is not smooth), but they belong to a more general class that
we are now going to define.

\begin{definition}
\label{almost}
A map $\Tr$ is said to be an {\it almost smoothing perturbation of the
  identity} if there exist smoothing functions $q^j\in\resto^i_l$ for
some $i,l\geq0$, and $S^k_1\in\Sc_1^k$, for some $k\geq 0$, such that
the following representation formula holds
\begin{equation}
\label{almost.1}
\Tr(\phi)=e^{q^j(\phi)JA_j}(\phi+S_1^k(\phi))\ .
\end{equation} 
\end{definition} 
\begin{remark}
\label{s.0}
The range of the smoothing map $S^k_1$ of equation \eqref{almost.1} is
not contained in $\V^0$, but in $\Hc^0$.
\end{remark}

\begin{remark}
\label{rem.sin}
An almost smoothing perturbation of the identity $\Tr$ is also almost smooth, namely
$\Tr\in\ac(\V,\V)$.  
\end{remark}

\begin{remark}
\label{s.3s}
Almost smoothing perturbations of the identity form a group of
continuous mappings, furthermore for any $j$ one has
$\Ph\relax_j(\Tr(\phi))=\Ph\relax_j(\phi)+R_2^1(\phi)$.
\end{remark}

\begin{proposition}
\label{pr.almost}
Let $H$ be a Hamiltonian function defined on $\Hc^{k_0}$ for some $k_0$;
assume that it is invariant under the symmetry group, namely that
$H(u)=H(e^{q^jJA_j}u)$, and consider $H_r:=i^*H$; let $\Tr$ be an
almost smoothing perturbation of the identity with $q^j\in\resto^1_2$ and $k=1$ (vanishing index
of the map $S$), then one has
\begin{equation}
\label{s.e.1}
H_r(\Tr(\phi))=H(\eta_{p_0-N}+S^1_2+\Pi_{p_0-N}( \phi+ S^1_1))\ ,
\end{equation}
with suitable maps $S^1_2$ and $S^1_1$.
\end{proposition}
\proof  First, remark that
for any choice of the scalar quantities $q^j$, and $p_j$,one has
\begin{equation}
\label{not.1.2}
\norma{(\pip e^{q^jJA_j}-e^{q^jJA\relax_j}\pip)\phi}_{\Hc^k}\leq
C |q|\norma{\phi}_{\Hc^{-l}} \ ,
\end{equation}
so that 
\begin{equation}
\label{e.cor.2}
\pip e^{q^jJA\relax_j}\phi=e^{q^jJA\relax_j}( \pip\phi+S^1_3)
\end{equation}

Write $p'=p\circ \Tr$, $N'=\Ph(\Tr(\phi))=N +R_2^1$, then one has
$p'=p_0-N'+R_2^1\circ\Tr= p_0-N'+R^1_2$. Consider now 
\begin{eqnarray*}
\eta_{p'}+\Pi_{p'}(\e^{q^jJA_j} (\phi+S^1_1))=
\eta_{p_0-N}+S_2^1+\Pi_{p_0-N}(\e^{q^jJA_j} (\phi+S^1_1))+S^1_2 
\\
=
\eta_{p_0-N}+S_2^1+\e^{q^jJA_j}\Pi_{p_0-N}
(\phi+S^1_1)+S^1_3
\\
=\e^{q^jJA_j}\left(\e^{-q^jJA_j}
\eta_{p_0-N}+\e^{-q^jJA_j}S_2^1+\Pi_{p_0-N}
(\phi+S^1_1)) \right)
\\
=\e^{q^jJA_j}\left(
\eta_{p_0-N}+S^1_2 +S_2^1+\Pi_{p_0-N}
(\phi+S^1_1)) \right)\ .
\end{eqnarray*}
Inserting in $H$ and exploiting its invariance under the group action
$\e^{q^jJA_j}$ one gets the result.\qed

\begin{lemma}
\label{l.6}
Let $s^l\in \resto^a_j$, $X\in\Sc^a_i $ $j\geq i\geq 1$, $a\geq 0$ be
smoothing functions, and consider the equation
\begin{equation}
\label{l.5.12}
\dot \phi= s^l(N,\phi)\Pi_{p_0}JA_l\phi+X(N,\phi)\ .
\end{equation}
Then for $|t|\leq 1$, the corresponding flow $\F_t$ exists in a
sufficiently small neighborhood of the origin of $\V^{d_A/2}$, and
for any $|t|\leq 1$ it is an almost smoothing perturbation of the
identity of the form
\begin{equation}
\label{flt}
\F_t=\e^{q^l(N,\phi,t)JA_l}(\phi+S(N,\phi,t))\ ,
\end{equation}
with $q^l(t)\in\resto^a_j$ and $S(t)\in\Sc^a_{i}$. Furthermore one
has
\begin{equation}
\label{flt1}
N(t)=N+S^a_{i+1}\ .
\end{equation}
\end{lemma}
\begin{remark}
\label{tei.r}
In sect. \ref{Normalform} we will study the case where $X$ is a
homogeneous polynomial of degree $r-1$ in $\phi,$ $r\geq 3$ and $s^l$
are homogeneous polynomials of degree $r$, then 
$\F_1\in\ac(\V,\V)$ can be expanded in Taylor series up to any finite
order and the remainder can be estimated. In particular, for any $k\geq 0$
there exists $\bar k$ s.t.
\begin{equation}
\label{sti.ta}
\norma{\F_1(\phi)-\left(\phi+X(\Ph(\phi),\phi)\right)}_{\Hc^k}\leq
C_k\norma{\phi}_{\Hc^{\bar k}}^r\ .
\end{equation}\end{remark}

\proof First rewrite \eqref{l.5.12} as
\begin{equation}
\label{flow.1}
\dot \phi=s^l(N,\phi)JA_l\phi+X+S^a_{j+1}=s^l(N,\phi)JA_l\phi+S^a_{i}\ ,
\end{equation}
then we rewrite the equation \eqref{flow.1}
in a more convenient way, namely we add a separate equation for the
evolution of $N$ and then we use a variant of Duhamel principle in
order to solve the system. Write
\begin{align}
\label{l.6.5}
\dot \phi =Y(N,\phi)\ ,\quad \dot N=Y_N(N,\phi)\ ,
\end{align}
where 
\begin{align*}
Y(N,\phi)&:=s^l(N,\phi)JA_l\phi+S^a_{i}
\\
Y_{Nk}&:=\langle\dot \phi;A\relax_k\phi\rangle= s_l\langle JA_
l\phi;A\relax_k\phi\rangle+ \langle S^a_{i};\phi\rangle = \langle
S^a_{i};\phi\rangle \ .
\end{align*}
In order to solve the system \eqref{l.6.5} we make the Ansatz $
\phi=e^{q^jJA_j}\psi $ with $q^j$ auxiliary variables. Compute $\dot
\phi$, and impose $\dot q^l=s^l(N,\phi)$, thus we get the system
\begin{align}
\label{sys.to.1}
\dot \psi=\e^{-q^jJA_j}S^a_{i}(N,\e^{q^jJA_j}\psi)\ ,
\ \dot q^l=s^l(N,\e^{q^jJA_j}\psi)\ ,
\ \dot N=Y_N(N,\e^{q^jJA_j}\psi)\ ,
\end{align} 
which is equivalent to \eqref{flow.1}. Fix now $k$ and consider such a
system in $\Hc^k\oplus\R^n\oplus\R^n$. Since the r.h.s. is smoothing,
for any $r$ there exists a neighborhood of the origin $\U\subset
\Hc^k\oplus\R^n\oplus\R^n$ in which it is of class
$C^r(\U;\Hc^k\oplus\R^n\oplus\R^n)$. It follows that it generates a
local flow which is also of class
$C^r((-\epsilon,\epsilon)\times\tilde \U;\Hc^k\oplus\R^n\oplus\R^n)$
(the first interval being that to which time belongs and $\tilde
\U\subset \U$). To show that in a small neighborhood of the origin the
flow is defined up time $1$ and to get the estimates ensuring the
structure \ref{almost.1}, just remark that we have
\begin{eqnarray*}
|\dot q |\leq C_l |N|^a\norma \psi_{\Hc^{-l}}^j
\ ,\ 
|\dot N|\leq C_l |N|^{a}\norma \psi_{\Hc^{-l}}^{i+1}
\ ,\ 
\norma{\dot \psi}_{\Hc^k}\leq C_{kl} |N|^a\norma \psi_{\Hc^{-l}}^i \ ,
\end{eqnarray*}
then the standard theory of a priori estimates of differential
equations gives the result. 
\qed
\begin{remark}
\label{smo.flow}
As we will see, in Darboux coordinates, the Hamiltonian vector field
of a smoothing Hamiltonian has the structure \eqref{l.5.12}, thus such
Hamiltonians generate a flow of almost smoothing perturbations of the
identity.
\end{remark}
\begin{remark}
\label{time.dep}
The result of Lemma \ref{l.6} immediately extends to the case where
$s^l$ and $X$ also depend on smoothly on time. Precisely
$X=X(t,N,\phi)$ must be such that $\forall r,l,k\geq 0$ there exists
an open neighborhood of the origin $\U_{rlk}\subset \R^n\oplus\V^{-k}$
s.t. $X\in C^r\left([-1,1]\times\U_{rlk},\V^l\right)$, and $s^l$ must
fulfill a similar property. 
\end{remark}
\begin{theorem}
\label{darboux}
(Darboux theorem) There exists an almost smoothing perturbation of the
identity 
\begin{equation}
\label{t.1.1}
\phi=\F(\phi')=e^{q^jJA_j}(\phi'+S^1_1(\phi'))
\end{equation}
with $q^j\in\resto^1_2$, such that $\F^*\Omega=\Omega_0$, i.e., in the
coordinates $\phi'$ one has
\begin{equation}
\label{t.1.5}
\Omega(\Phi'_1;\Phi'_2)=\langle E\Phi'_1;\Phi'_2\rangle\ .
\end{equation}
Correspondingly the Hamilton equations of a Hamilton function
$H(\phi')$ have the form
\begin{equation}
\label{t.1.6}
\dot \phi'=\Pi_{p_0}J\nabla H(\phi')
\end{equation}
\end{theorem}

The rest of the section is devoted to the proof Theorem \ref{darboux}.

We recall that in standard Darboux theorem the transformation
introducing canonical coordinates is constructed as follows.  Denote
$\Omega_0:=\Omega_{\phi}\big|_{\phi=0}$,
$\tilde\Omega:=\Omega-\Omega_0$ and $\Omega_t:=\Omega_0+t\tilde\Omega$.
Let $\alpha$ be a 1-form such that $\tilde\Omega=\di\alpha$ and let
$Y_t$ be such that $ \Omega_t(Y_t,.)=-\alpha$. Let $\F_t$ be the
evolution operator of $Y_t$ (we will prove that it is formed by almost
smoothing perturbations of the identity), then
\begin{equation}
\label{i.9a}
\frac{d}{dt}\F^*_{t}\Omega_t=\F^*_{t}(\lie_{Y_{t}}\Omega_t)+ \F^*_{t}\frac 
d{dt} \Omega_t=\F^*_{t}(-d\alpha+\tilde\Omega)=0\ ,
\end{equation}
so $\F^*_1\Omega\equiv \F^*_1\Omega_1=\Omega_0$, and $\F:=\F_1$ is the
wanted change of variables. We follow such a scheme, by adding the
explicit estimates showing that all the objects are well defined.

\begin{remark}
\label{darb.rem}
It will turn out that $\F_t:\V^k\to\V^k$ is differentiable with
respect to time at the points $\phi\in\V^{k_0}$, with a sufficiently
large $k_0$. As a consequence one immediately see that the Darboux
procedure is valid for solutions with initial data in $\V^{k_0}$. The
extension to general initial data is obtained by density.
\end{remark}

First we compute the expression of the symplectic form in the
coordinates introduced by lemma \ref{impl.1}. In order to
simplify the computation we will first compute $\Omega^0:=i_0^*\omega$
with $i_0$ the map \eqref{i.7a}. It is also useful to compute a 1-form $\Theta^0$ such that $\di
\Theta^0=\Omega^0$. Subsequently we compute $\Omega=i^*\omega$ and a
potential 1-form for $\Omega$ by inserting the expression of
$p=p(N,\phi)$.

\begin{lemma}
\label{l.1}
Define the 1-form $\Theta^0$ by
\begin{equation}
\label{l.1.1}
\Theta^0(P,\Phi)= \frac{1}{2}\langle E\Pi_p\phi;\dep{\Pi_p}j\phi
\rangle P_{j}+\langle E\Pi_p\phi;\Phi\rangle\ 
\end{equation}
(by this notation we mean that the r.h.s. gives the action of the form
$\Theta^0$ at the point $(p,\phi)$ on a vector $(P,\Phi)$), 
then one has $\di \Theta^0=\Omega^0\equiv i^*_0\omega$, and
therefore
\begin{align}
\label{l.1.2}
\Omega^0\left((P_1,\Phi_1);(P_2,\Phi_2) \right)= \left\langle E\depp
j\phi,\depp i\phi\right\rangle P_{1j}P_{2i} +\left\langle E\Pi_p\depp
j\phi ; \Phi_2 \right\rangle P_{1j}
\\
\nonumber
- \left\langle E\Pi_p\depp
j\phi ; \Phi_1 \right\rangle P_{2j}+\langle E\pip \Phi_1;\Phi_2
\rangle\ .
\end{align}
\end{lemma}
\proof We compute $i_0^*\theta$, where $\theta=\langle Eu;.\rangle/2$
is such that $\omega=\di \theta$.  By writing $u=i_0^*(p,\phi)$, one
has
\begin{align}
\label{i.25}\dep uj=\dep{\eta_p}j+\dep{\Pi_p}j\phi \ ,\quad
(\di_\phi i_0)\Phi=\Pi_p\Phi
\end{align}
so, taking $\theta=\frac{1}{2}\langle Eu;.\rangle$, one has
$$
(i^*_0\theta)(P,\Phi)=\frac{1}{2}\left\langle Eu;\dep uj
\right\rangle P\relax_j+\frac{1}{2}\left\langle Eu;\di_\phi
i_0\Phi\right\rangle \ .
$$
We compute the first term, which coincides with
\begin{align}
\label{i.25b}
2\theta\left(\dep uj\right)
&=\langle
E(\eta_p+\Pi_p\phi); \dep{\eta_p}j+\dep{\Pi_p}j\phi\rangle 
\\ \nonumber
&=\langle E\eta_p;\dep {\eta_p}j\rangle+\langle
E\eta_p;\dep{\Pi_p}j\phi\rangle +\langle E\Pi_p\phi;{\dep
  {\eta_p}j}\rangle +\langle E\Pi_p\phi;\dep{\Pi_p}j\phi\rangle\ .
\end{align}
Now, the third term vanishes due to the definition of
$\pip$. Concerning the first term, there exists a function $f_0$ such
that $\dep{f_0}j= \langle E\eta_p;\dep {\eta_p}j\rangle$, indeed, from
the isotropy property (B3) one has
$$ \dep\null j\langle E\eta_p;\dep {\eta_p}i\rangle=\dep\null i\langle
E\eta_p;\dep {\eta_p}j\rangle\ .
$$ Finally, defining $\displaystyle{f_1(p,\phi)=\langle E\eta_p;\pip
  \phi\rangle}$, the second term of \eqref{i.25b} turns out to be
given by $\dep{f_1}j$, so we have
$$
2\theta\left(\frac{\partial u}{\partial p_j}\right)=\langle
E\Pi_p\phi;\dep{\Pi_p}j\phi\rangle+\frac{\partial(f_0+f_1)}{\partial p_j}\ .
$$
We compute now $\langle Eu; (\di_\phi i_0)\Phi\rangle$. We have
$$ 2\langle Eu; (\di_\phi i_0)\Phi\rangle=\langle
E(\eta_p+\Pi_p\phi);\Pi_p\Phi\rangle= \langle
E\Pi_p\phi;\Pi_p\Phi\rangle+(\di_{\phi}f_1)\Phi\ ,
$$
from which $i_0^*\theta=\Theta^0+\di(f_0+f_1)$, and therefore
$\Omega^0=\di\Theta^0$. 

We compute now explicitly $\di\Theta^0$. Denote $ \Theta^0=
\Theta^{0j} \di p\relax_j+\langle\tilde \Theta^0_\phi;.  \rangle $,
then the computation of $\frac{\partial \Theta^{0i}}{\partial
  p_j}-\frac{\partial \Theta^{0j}}{\partial p_i}$ is trivial and gives
the term proportional to $P_{1j}P_{2i}$ in \eqref{l.1.2}. Also the
computation of the term containing $\Phi_1,\Phi_2$ is trivial and is
omitted. We come to the $P,\Phi$ terms. When applied to a vector
$\Phi$ it is given by
\begin{align}
\label{l.1.4}
&\langle \dep{ \Theta^0_\phi}j;\Phi\rangle-(\di_\phi
\Theta^0_j)\Phi 
\\
&= \frac{1}{2}\left( \left\langle E\depp j\phi;\Phi
\right\rangle -\left\langle E\pip\Phi; \depp j\phi
\right\rangle- \left\langle E\pip\phi;\depp j\Phi
\right\rangle\right)\ ,
\end{align}
which is the scalar product of $\Phi$ with a half of the vector
\begin{eqnarray*}
E\depp j\phi+\pip^* E\depp j\phi-\dep{\Pi_p^*}j E\pip \phi
= E\left(\depp j+\pip\depp j -\depp j\pip \right)\phi
\\
= E\left(\pip\depp j +\depp j\pip+ \pip\depp j -\depp j\pip\right)\phi
=2 E\pip \depp j\phi\ ,
\end{eqnarray*}
which immediately gives the thesis. \qed

\begin{lemma}
\label{l.2}
In the coordinates of lemma \ref{impl.1} the symplectic form
$\Omega=i^*\omega$ takes the form $\Omega(\Phi_1,\Phi_2)=\langle
O\Phi_1;\Phi_2\rangle$ with $O$ given by
\begin{align}
\label{l.2.1}
O\Phi= a^{ij}\langle\nabla p\relax_i;\Phi\rangle \nabla p\relax_j +\langle\nabla
p\relax_j;\Phi\rangle \Pi_{p_0}^*E\pip\depp j\phi
\\
- \left\langle E\pip\depp
j\phi;\Phi\right\rangle\nabla p\relax_j+\Pi_{p_0}^* E\pip\Phi  
\end{align}
where 
\begin{equation}
\label{l.2.2}
a^{ij}:=\langle E\depp i\phi;\dep{\Pi_p}j\phi
\rangle\ ,
\end{equation}
and $p=p(N,\phi)$.  Moreover one has $\Omega=\di \Theta$ with
\begin{equation}
\label{l.2.3}
\Theta(\Phi)=\frac{1}{2}\langle E\pip\phi;\Phi\rangle+ \frac{1}{2}
\left\langle E\pip\phi;\depp j\phi \right\rangle \langle\nabla
p\relax_j;\Phi\rangle
\end{equation}
\end{lemma}
\proof The expression of $\Omega$ and $\Theta$ are obtained by taking
\eqref{l.1.1} and  \eqref{l.1.2} and inserting the expression of
$p=p(\phi)\equiv p(\Ph(\phi),\phi)$ and substituting 
$
P_{1,2 j}=\langle\nabla p\relax_j;\Phi_{1,2}\rangle
$, thus
the thesis follows from a simple computation.\qed

\begin{remark}
\label{r.t}
One can define $\Omega_t=\langle O_t.;.\rangle$ and $\alpha=\langle
V;.\rangle$ with
\begin{align}
\label{r.t.1}
O_t&=E+t\left[\Pi_{p_0}^* E\pip-E+a^{ij}\langle\nabla
  p\relax_i;.\rangle \nabla p\relax_j\right. 
\\ 
\nonumber &+\left.\langle\nabla
  p\relax_i;.\rangle \Pi_{p_0}^*E\pip\depp j\phi\right.  \left.  -
  \left\langle E\pip\depp j;.\right\rangle\nabla p\relax_j \right]
\\ 
V&=-\frac{1}{2} \left(\Pi_{p_0}^*E\pip\phi -E\phi+\left\langle
E\pip\phi;\depp j\phi \right\rangle \nabla p\relax_j \right)
\end{align}
\end{remark}
In order to find the normalizing vector field $Y_t$, we have to solve the
equation
\begin{equation}
\label{nor.1}
O_t Y=-V
\end{equation}
where we omitted the index $t$ from $Y$. We start now the discussion
of such an equation.

First we have the lemma: 
\begin{lemma}
\label{l.3}
Define $D_t:=E+t(\Pi_{p_0}^*E\pip -E)$, then it is skew-symmetric;
furthermore, provided $|p-p_0|$ is small enough $\exists S_t$
fulfilling $S_t(\V^{-\infty})^*\subset\V^\infty$, s.t. $
D_t^{-1}=J+S_t$ and
\begin{equation}
\label{i.32a}
\norma{S_t\Phi}_{\Hc^k}\leq C_{k,l}\left|p-p_0\right|\norma{\Phi}_{-l}
\leq C'_{k,l}|N|\norma{\Phi}_{-l} \ .
\end{equation}
\end{lemma}
\proof First remark that since $D_t$ acts on $\V^k$, the term
$\Pi_{p_0}^*E\pip$ can be rewritten in the more symmetric form
$\Pi_{p_0}^*E\pip\Pi_{p_0} $, from which, using \eqref{r.2.1} one
immediately sees skew-symmetry.  We have now $D_t=E+t\tilde D$, where
$\tilde D:=\Pi_{p_0}^*E\Pi_p-E$, which is smoothing and fulfills an
inequality equal to \eqref{i.32a}. Then $D_t=E(\uno+tJ\tilde D)$, and
by Neumann formula one gets
$$
D_t^{-1}=J-tJ\tilde D\sum_{k\geq0}(-1)^k(tJ\tilde D)^kJ
$$
and the thesis. The second of \eqref{i.32a} follows from Lemma
\ref{impl.1}, item (2).\qed

\begin{lemma}
\label{l.4}
The solution of equation \eqref{nor.1} has the form
\begin{equation}
\label{l.4.1}
Y(p,\phi)=s^l(N,\phi,t)\Pi_{p_0}JA\relax_l\phi+S^1_1(t)
\end{equation}
where $s^j(t)\in\resto^1_2$, $S^1_1(t)\in\Sc^{1}_1$ and
$S^1_1$ are smoothly dependent on 
$t\in[0,1]$. 
\end{lemma}
\proof First write explicitly \eqref{nor.1}
introducing, for short, the notations
\begin{eqnarray*}
b\relax_i:=\langle\nabla p\relax_i,Y\rangle\ ,\quad W^i=\Pi_{p_0}^*E\pip\depp i\phi\ ,
\end{eqnarray*}
so that it takes the form
\begin{align}
\label{l.4.5}
D_tY+t\left(a^{ij}b\relax_i\nabla p\relax_j+b\relax_iW^i-\left\langle W^j,Y
\right\rangle \nabla p\relax_j  \right) 
\\
\nonumber
=\frac{1}{2}[E\pip
  -E]\phi-\frac{1}{2} \langle \phi; W^j\rangle\nabla p\relax_j\ .
\end{align}
Applying $D_t^{-1}$ and reordering the formula one gets
\begin{align}
\label{l.4.6}
Y=-t\left(a^{ij}b\relax_i\dm \nabla p\relax_j+b\relax_i\dm W^i-\langle
W_j;Y\rangle\dm\nabla p\relax_j\right)
\\
\nonumber
-\frac{1}{2} \langle \phi
;W^j\rangle\dm \nabla p\relax_j
+\frac{1}{2}\dm [E\pip -E]\phi\ .
\end{align}
Denote 
$
\gamma\relax_{ij}:=\left\langle\nabla p\relax_i;\dm \nabla
p\relax_j\right\rangle$, $\beta\relax_i^j:=\langle \nabla p\relax_i;\dm W^j\rangle
$
and remark that $\gamma_{ij}$ is smoothing, since it is given by
\begin{align*}
\left\langle\nabla p\relax_i;J \nabla
p\relax_j\right\rangle+\left\langle\nabla p\relax_i;S_t \nabla
p\relax_j\right\rangle=M\relax_i^l\langle A\relax_l\phi;JA\relax_k\phi\rangle M\relax_j^k+\left\langle\nabla p\relax_i;S_t \nabla
p\relax_j\right\rangle
\\
=\left\langle\nabla p\relax_i;S_t \nabla
p\relax_j\right\rangle  \  .
\end{align*} 
Also $\beta\relax_i^j$ is clearly smoothing. Now, taking the scalar
product of \eqref{l.4.6} with $\nabla p_i$ and $W^l$ respectively, one
has
\begin{eqnarray*}
b\relax_i&=&-t\left( \gamma\relax_{ij}a^{ik}b\relax_k+\beta\relax_i^j b\relax_j-\langle
W^j;Y\rangle\gamma\relax_{ij} \right) \\&-& \gamma\relax_{ij}\frac{1}{2}\langle
\phi;W^j\rangle+\frac{1}{2}\langle \nabla p\relax_i;\dm
\left[E\pip-E\right]\phi\rangle
\\
\langle W^l;Y\rangle&=& -t\left(a^{ij}b\relax_i\beta\relax_j^l+b\relax_i\langle W^l;\dm
W^i\rangle- \langle W^j; Y\rangle\beta\relax_j^l\right)
\\
&-&\frac{1}{2}\langle\phi;W^j\rangle\beta\relax_j^l+\frac{1}{2}\langle W^l;\dm
[ E\pip-E]\phi\rangle\ ,
\end{eqnarray*}
which is a linear system for $b\relax_i$ and $\langle W^l;Y\rangle$. Solving
it and 
analyzing the solutions one gets
\begin{align}
\label{l.4.10}
b\relax_i&= \frac{1}{2}\langle \nabla p\relax_i;\dm
\left[E\pip-E\right]\phi\rangle+{\rm h.o.t.}\ ,
\\
\label{l.4.11}
\langle W^l;Y\rangle&=\frac{1}{2}\langle W^l;\dm
[ E\pip-E]\phi\rangle+{\rm h.o.t.}\ ,
\end{align}
where h.o.t are also regularizing. Substituting in
\eqref{l.4.6} one gets a formula for $Y$. Then one has that such an
$Y$ actually fulfills \eqref{l.4.5}, and is thus the wanted solution
of \eqref{nor.1}. 
From these formulae one immediately has
\begin{align}
\label{l.4.12}
| b\relax_i|\leq C\left|p-p_0\right|\norma{\phi}_{\Hc^{-k}}^2 \leq
CN \norma{\phi}_{\Hc^{-k}}^2
\\
\label{l.4.13}
|\langle W^l;Y\rangle  |\leq
C\left|p-p_0\right|^2\norma{\phi}_{\Hc^{-k}}^2
\leq
CN^2 \norma{\phi}_{\Hc^{-k}}^2
\end{align}
where we used item (2) of Lemma \ref{impl.1}.
 To get the formula \eqref{l.4.1} and the corresponding estimates,
 define the function $s^l$ to be the coefficient of $\dm\nabla p\relax_l$ in
 \eqref{l.4.6}, and remark that $\dm\nabla
 p\relax_l=\Pi_{p_0}JA\relax_l\phi$+smoothing terms. Then it is easy
 to conclude the proof.\qed

Theorem \ref{darboux} is an immediate consequence of Lemma \ref{l.6}. 

\section{ The Hamiltonian in Darboux coordinates}\label{hami}

Concerning the smoothness and the structure of the nonlinear part of
the Hamiltonian we make the following assumption

\begin{itemize}
\item[(P2)] 
The map 
\begin{equation}
\label{H2.1}
\eta\mapsto \nabla H_P(\eta)\ 
\end{equation}
 is of class $\ac(\Hc,\Hc)$. 

Denote 
$$
\di^2 H_P(\eta)(\Phi;
\Psi)\equiv\frac{1}{2}\langle B(\eta)\Phi;\Psi\rangle\ ,
$$ then the map $(\eta,\Phi)\mapsto B(\eta)\Phi$ is of class
$\ac(\Kcs;\Hc)$, $\Kcs=\left\{\Hc^k\oplus\Hc^k\right\}_{k\in\Z}$.
\end{itemize}

\subsection{The Hamiltonian in Darboux coordinates}\label{tay}

We  introduce the coordinates of the Darboux theorem \ref{darboux}.  To this end
we exploit proposition \ref{pr.almost} from which one gets:

\begin{proposition}
\label{p.ham.1}
In the Darboux coordinates introduced by theorem \ref{darboux} the
Hamiltonian $H_r\circ\F$ has the form
\begin{align}
\label{e.ham.1}
H_r(\F(\phi))&=H_L+H_N\ ,
\\
H_L&:= \Phzero (\phi)+\frac{1}{2 }\di^2
H_P(\eta_{p_0-N})(\phi,\phi) -\lambda^j(p_0)\Ph_j+D(N)+(R^1_2)_{lin}
\\
\nonumber
H_N&:=R^1_3+{H^3_P}(\eta_{p_0-N},\phi+S^1_1)
\end{align}
where $(R^1_2)_{lin}$ is a smoothing quadratic polynomial in $\phi$
and  
\begin{align}
\label{e.ham3}
D(N)&:=H(\eta_{p_0-N})-[H(\eta_{p_0})- \frac{\partial H
  }{\partial p_j}(\eta_{p_0})N_j ]\ ,
\\
\label{e.ham.8.1}
H^3_P(\eta,\phi)&:=H_P(\eta+\phi)-[H_P(\eta)+\di
  H_P(\eta)\phi+\frac{1}{2}\di^2H_P(\eta)(\phi,\phi)   ] \ .
\end{align}
\end{proposition}
\proof Exploiting proposition \ref{pr.almost} one has to study
\begin{eqnarray*}
\Ph_0(\eta_{p_0-N}+S^1_2+\Pi_{p_0-N}(\phi+S^1_1))=
\Ph_0(\eta_{p_0-N}+\Pi_{p_0-N}(\phi+S^1_1))+R^1_2 
\\ 
= \Ph_0(\eta_{p_0-N})+\left\langle A_0\eta_{p_0-N};
\Pi_{p_0-N}(\phi+S^1_1) \right\rangle+\Phzero(\Pi_{p_0-N}(\phi+S^1_1)
)+R^1_2\ ,
\end{eqnarray*}
and
\begin{eqnarray*}
H_P(\eta_{p_0-N}+S^1_2+\Pi_{p_0-N}(\phi+S^1_1))
\\
=H_P(\eta_{p_0-N})+\di
H_P(\eta_{p_0-N}) \left(S^1_2+\Pi_{p_0-N}(\phi+S^1_1) \right) 
\\
+\frac{1}{2}\di^2
H_P(\eta_{p_0-N}) \left(S^1_2+\Pi_{p_0-N}(\phi+S^1_1) \right)
^{\otimes 2}
\\
 +H^3_P(\eta_{p_0-N};S^1_2+\Pi_{p_0-N}(\phi+S^1_1) )\ .
\end{eqnarray*}
Consider first the terms linear in $\phi$: they are given by
\begin{eqnarray*}
\left\langle A_0\eta_{p_0-N};
\Pi_{p_0-N}(\phi+S^1_1) \right\rangle + \di
H_P(\eta_{p_0-N}) \left(\Pi_{p_0-N}(\phi+S^1_1) \right)
\\
=\lambda^j(p_0-N)\left\langle A_j\eta_{p_0-N};\Pi_{p_0-N}(\phi+S^1_1)
\right\rangle =0
\end{eqnarray*}
where we exploited eq. \eqref{e.ham.4}. Thus we have
\begin{eqnarray*}
H_r(\F(\phi))= H(\eta_{p_0-N})+\Phzero (\phi+S^1_1)+\frac{1}{2 }\di^2
H_P(\eta_{p_0-N})(  \phi+S^1_1,\phi+S^1_1) 
\\
+R^1_2+{H^3_P}(\eta_{p_0-N},\phi+S^1_1)
\\
= H(\eta_{p_0-N})+\Phzero (\phi)+R^1_2+\frac{1}{2 }\di^2
H_P(\eta_{p_0-N})(  \phi,\phi) 
\\
+R^1_2+{H^3_P}(\eta_{p_0-N},\phi+S^1_1)\ .
\end{eqnarray*}
Finally we have to rewrite in a suitable form the function
$H(\eta_{p_0-N})$. Denote $f(p):=H(\eta_p)$. Expanding at $p_0$ and
using equation
\eqref{e.ham.4}, one has
\begin{equation}
\label{e.ham.9}
\dep fj(p_0)= \di H(\eta_{p_0})
\dep{\eta_{p_0}}j
=\lambda^k\langle
  A\relax_k\eta_{p_0};
  \dep{\eta_{p_0}}j\rangle
=\lambda^k\delta_k^j=\lambda^j\ ,
\end{equation}
from which the thesis immediately follows.
\qed

\begin{remark}
\label{vf.3}
Define $X_P:=J\nabla H_P$, 
and, for fixed $\eta\in \Hc^\infty$, 
\begin{equation}
\label{vf.31}
X_P^2(\eta;\phi):=X_P(\eta+\phi)-[X_P(\eta)+\di
  X_P(\eta)\phi]\ ,
\end{equation}
 then one has $J\nabla
H_P^3(\eta;\phi)=X_P^2(\eta,\phi)$. This can be seen by writing the
definition of Hamiltonian vector field.  
\end{remark}
\begin{remark}
\label{field}
The  Hamilton vector field of {\bf $H_r\circ\F$, which from now on will be
simply denoted by $H$} is given by 
\begin{align}
\label{ham.eq}
\dot \phi=\Pi_{p_0}
J[A_0\phi+V_N\phi+S^1_1(N,\phi)+X^2_P(\eta_{p_0-N},\phi)
  -\lambda^j(p_0)A_j\phi] 
\\
\nonumber
+w^j(N,\phi)\Pi_{p_0} JA_j\phi\ ,
\end{align}
where
\begin{equation}
\label{ham.eq.1}
w^j:=\lambda^j(p_0)-\lambda^j(p_0-N)+\frac{1}{2}\langle\frac{\partial
  V_N}{\partial N_j}\phi;\phi \rangle +S^0_2+\frac{\partial
  H_N}{\partial N_j}
\end{equation}
and $V_N$ is the operator such that 
\begin{equation}
\label{vN}
d^2H_P(\eta_{p_0-N})(\phi,\phi)=\frac{1}{2}\langle V_N\phi;\phi \rangle\ ,
\end{equation}
so that $V_N\phi=\di X_P(\eta_{p_0-N})\phi$.
\end{remark}

\subsection{Adapted
  coordinates}\label{linear} 

Consider the quadratic part of the original Hamiltonian at
$\eta_{p_0}$, namely
\begin{equation}
\label{e.lin.233}
H_{L0}(u):=\Phzero(u)+\frac{1}{2}\di^2H_P(\eta_{p_0})(u,u)-\lambda^j(p_0)\Ph_j(u)
\end{equation}  
 Denote 
$B:=A_0+V_0-\lambda^j(p_0)A_j$, $L_0:=JB$,
so that $H_{L0}(u)=\langle u;Bu\rangle/2=\langle EL_0 u;u\rangle/2 $.
Making the Ansatz
$u=\e^{t\lambda^jJA_j}(\eta_{p_0}+\chi)$ and linearizing in $\chi$ the
Hamilton equations of \eqref{e.ham}, one gets that $\chi$ satisfies $\dot
\chi=L_0\chi$.
\begin{lemma}
\label{lem.ker}
The generalized kernel of $L_0^*$ contains the vectors
\begin{equation}
\label{e.lin.333}
A_j\eta_{p_0}\ ,\quad E\frac{\partial \eta_{p_0}}{\partial
  \lambda^j}\ .
\end{equation}
\end{lemma}
\proof First, one immediately sees that $JA_j\eta_{p_0}\in Ker (L_0)$,
then exploiting the equation \eqref{e.ham.4} for the ground state one
sees that $B \frac{\partial \eta_{p}}{\partial \lambda^j}=A_j\eta_p$,
which implies $[L_0^*]^2E\frac{\partial \eta_{p_0}}{\partial
  \lambda^j}=[(JB)^*]^2E\frac{\partial \eta_{p_0}}{\partial \lambda^j}
=0 $, from which one immediately sees that the generalized kernel of
$(L_0)^*$ contains the vectors \eqref{e.lin.333}. \qed

We assume

\begin{itemize}
\item[(L1)] The generalized Kernel of $L_0^*$ is
  $2n$-dimensional. Furthermore $L_0\big|_{\V^j}$ is an isomorphism
  between $\V^j$ and $\V^{j-d_0}$.
\item[(L2)] $\langle B\phi,\phi\rangle >0$, $\forall
  \phi\not=0$, $\phi\in\V^{d_0/2}$.  
\item[(L3)] the essential spectrum of $L_0\big|_{\V^0}$ is
  $\bigcup_{\pm}\pm \im [{\Omega},+\infty)$. The rest of the spectrum
  consists of purely imaginary eigenvalues $\pm\im \omega_j$, that we
  order as follows $0<\omega_1\leq\omega_2\leq...\leq
  \omega_K$. Furthermore the corresponding eigenfunctions $v_{j\pm}$
  are smooth, namely $v_{j\pm}\in(\V^{\infty})^{\otimes\C}$.
\end{itemize}

In order to perform the dispersive estimates we will also have to
avoid boundary resonances. Let $r_t$
be the smallest integer such that $r_t\omega_1\geq\Omega$. We assume
that
\begin{itemize}
\item[(L4)] One has $\omega\cdot k\not=\Omega$, $\forall k\in\Z^K$
  s.t. $|k|\leq 2r_t$.
\end{itemize}

We normalize the eigenfunctions in such a way that 
\begin{align}
\label{vjort}
\langle Ev_{j\pm},v_{k\pm}\rangle=0\ ,\ \forall j,k\ ,\quad \langle
Ev_{j+},v_{k-}\rangle=-\im \delta_{jk}\ , \quad
\overline{v_{j+}}=v_{j-}\ , 
\end{align}
which is always possible since \eqref{vjort} are the standard ``symplectic
orthogonality'' relations of the eigenfunctions of the operator $L_0$
(which is skew with respect to the symplectic form). In particular the
second of \eqref{vjort} is a consequence of (L2).

We now introduce coordinates $(\xi_j,\phi_c)$  by
\begin{equation}
\label{cooc}
\phi=\sum_{j=1}^{K}(\bar\xi_jv_{j+}+ \xi_jv_{j-})+\phi_c\ ,
\end{equation}
where $\phi_c$ is such that $\langle Ev_{j\pm};\phi_c\rangle=0$
$\forall j$'s.  Explicitly one has 
\begin{align}
\label{exi}
\xi_j=\im\langle Ev_{j+},\phi\rangle\ ,\quad \bar \xi_j=- \im\langle
Ev_{j-},\phi\rangle\ ,
\\
\label{exi2}
\phi_c:=P_c\phi:=\phi-\sum_{j,\pm}\pm \im \langle
Ev_{j\pm};\phi\rangle v_{j\mp}\ .
\end{align}
In these coordinates the phase space becomes 
\begin{equation}
\label{phaW}
(\xi,\phi_c)\in\C^K\oplus \W^j\ ,\quad \W^j:=P_c\V^j\ .
\end{equation}

As usual it is often useful to consider the
variables $\bar\xi_j$ as independent from the $\xi_j$'s.  Often we
will also denote $
\phi_d:=\sum_{j=1}^{K}(\bar\xi_jv_{j+}+
\xi_jv_{j-})$.
In these coordinates the Hamiltonian vector field of a Hamiltonian
function $H$ takes the form
\begin{align*}
\dot \xi_j=-\im \frac{\partial H}{\partial \bar\xi_j}\ ,\quad 
\dot \phi_c=J\nabla_{\phi_c}H\ ,
\end{align*}
and the main term of the quadratic part of the Hamiltonian
\eqref{e.ham} takes the form
\begin{equation}
\label{e.lin.12}
H_{L0}=\sum_{l}\omega_l|\xi_l|^2+\frac{1}{2}\langle EL_c\phi_c;\phi_c\rangle\ ,
\end{equation}
where $L_c:=P_cL_0P_c$. Concerning the momenta $\Ph_j$ one has
\begin{equation}
\label{3s}
\Ph\relax_j(\xi,\phi_c)=\frac{1}{2}\langle \phi_c ;
\A\relax_j\phi_c\rangle+P_1^j(\phi_c,\xi)+P_2^j(\xi)\ ,
\end{equation}
where
\begin{equation}
P_1^j(\phi_c,\xi):= \sum_{|\alpha|+|\beta|=1} \langle \phi_c;
E\Phi^j_{\alpha\beta}  \rangle\xi^\alpha\bar
\xi^\beta
\ ,\quad 
P_2^j(\xi):=\sum_{|\alpha|+|\beta|=2} A_{\alpha\beta}^j \xi^\alpha\bar
\xi^\beta\ ,
\end{equation}
and $\Phi_{\alpha\beta}^j\in(\W^\infty)^{\otimes \C}$,
$A_{\alpha\beta}^j$ are suitable functions and complex numbers, while
$\A_j:=P_cA_jP_c$. By a small abuse of notation, in the following we
will always denote in the same way $\W^j$ and its complexification. We
will also identify the two scales $\V\equiv\left\{\V^j\right\}$ and
$\left\{\R^K\oplus \W^j\right\}$. 
Denote by $M\relax_j$ the function
\begin{equation}
\label{eq.2}
M\relax_j:=\frac{1}{2}\langle\phi_c ;
 \A\relax_j\phi_c\rangle \ .
\end{equation}

Since $P_1$ and $P_2$ are smoothing, in the following the quantities
$M_j$ will play the role that in the previous sections was plaid by
the quantities $N_j$.

{\bf In the following we will substitute the classes $\resto^k_l$ by
  similar classes in which the functions $N$ are substituted by the
  functions $M$ and similarly for the classes $\Sc^k_l$}.

In order to make the translation we remark that, if $F\in \resto^1_2$
{\bf old classes}, then we have
\begin{equation}
\label{new.4}
F=R^1_2+R^0_4 \quad {\bf \text{new\ classes}}
\end{equation}
as it immediately follows from \eqref{3s}. Similarly one has
\begin{equation*}
\text{old}\quad  S^1_1\null `` = ''  S^1_1+S^0_2\ ,\quad
\text{old}\quad  R^1_2\null ``=''R^1_2+R^0_3\ .
\end{equation*}

\begin{remark}
\label{commu}
From the definition of the operators $\A\relax_j$ we have that they do not
fulfill assumption (S1) and (S2), instead they fulfill
\begin{equation}
\label{com.scrpt}
[\A\relax_j;\A\relax_k]=S\relax_{jk}\ ,\quad [\A\relax_k;J]=S\relax_k
\end{equation}
with $S\relax_{jk} $ and $S\relax_k$ smoothing operators. 
\end{remark}

\begin{lemma}
\label{lemlin.1}
One has
$ \e^{tJ\A_j}\phi_c=e^{tJA_j}\phi_c+S(t)\phi_c$,
where $S$ is a smoothing family of operators which fulfills
$
\norma{S(t)\phi_c}_k\sleq |t|\norma{\phi_c}_{-l}
$.
\end{lemma}
\proof Write explicitly the equation $\dot \phi=JA_j\phi$ (which
defines $\e^{tJA}$):
\begin{align}
\label{eq.lin}
\dot\phi_c &= J\A\relax_j\phi_c +\sum_{|\alpha+\beta|=1} \xi^\alpha\bar
\xi^\beta \Phi_{\alpha\beta}^j\ ,
\\
\dot \xi_k &= -\im\left(
\sum_{|\alpha+\beta|=2}\frac{\beta_k }{\bar\xi_k}
A^j_{\alpha\beta}\xi^\alpha\bar\xi^\beta+
\sum_{|\alpha+\beta|=1}\frac{\beta_k }{\bar\xi_k} \left\langle
E\Phi^j_{\alpha\beta};\phi_c \right\rangle  \xi^\alpha\bar\xi^\beta\right)\ .
\end{align}
Using Duhamel formula (with $J\A_j$ as principal part), one gets the
thesis.
\qed

With the new notations and classes one has that the Hamiltonian of the
system takes  the form
\begin{align}
\label{H.new}
H&=H_L+H_N \ ,\quad H_L=H_{L0}+H_{L1}+D(M)\ ,
\\
\label{H.new4}
H_{L0}&:=\Ph_0(\phi)+\frac{1}{2}\langle V_0\phi;\phi
\rangle-\lambda^j(p_0)
\Ph_j(\phi)=\sum_{l}\omega_l|\xi_l|^2+\frac{1}{2}\left\langle
EL_c\phi_c;\phi_c\right\rangle \ ,
\\
\label{H.new1}
H_{L1}&=\frac{1}{2}\left\langle (V_M-V_0)\phi;\phi \right\rangle
+(\resto^1_2) _{lin}
\\
\label{H.new2}
H_N&:=R^1_3+R^0_4+H^3_P(\eta_{p_0-M};\phi+S^1_1+S^0_2)\ ,
\end{align}
where $V_M$ is the operator $V_N$ evaluated at $N=M$.

Of course, this is also true for $H\circ \Tr$ with any $\Tr$ almost
smoothing perturbation of the identity of the form \eqref{t.1.1}.

In the following we will denote by $X_N$ the vector field of $H_N$
computed at constant $M$, i.e. as if $M$ were independent of
$\phi_c$.

\section{Normal Form}\label{Normalform}

First we define what we mean by normal form. 

\begin{definition}
\label{def.nor.for}
A function $Z(M,\xi,\phi_c)$, of class $Z\in \ac(\Kcs,\R)$,
$\Kcs:=\R^n\oplus \V$, will be said to be in normal form at order
$r$, if the following holds
\begin{align}
\label{eq.nor.for.1}
\left\{ \omega\cdot(\mu-\nu)\not=0\ \&\ |\mu|+|\nu|\leq
r\right\}\Longrightarrow \frac{\partial^{|\mu|+|\nu|
  }Z}{\partial\xi^{\mu}\partial\bar \xi^\nu}(M,0)=0
\\
\label{eq.nor.for.2}
\left\{ |\omega\cdot(\mu-\nu)|<\Omega\ \&\ |\mu|+|\nu|\leq
r-1\right\}\Longrightarrow \di_{\phi_c}\frac{\partial^{|\mu|+|\nu|
  }Z}{\partial\xi^{\mu}\partial\bar \xi^\nu}(M,0)=0
\end{align}
The derivatives with respect to $\phi_c$ have to be computed at
constant $M$, i.e. as if $M$ were independent of such quantities. 
\end{definition}
\begin{theorem}
\label{main.nf}
For any $r\geq 2$ there exists an almost smoothing, canonical
perturbation of the identity 
$\Tr_r(\phi)=\e^{q^lJA_l}(\phi+S(\phi))$, with $q^l\in
\resto^0_{2}$, and $S\in \Sc^1_1\cup\Sc^0_2$, such that
$H\circ\Tr_r$ is in normal form at order $r$.
\end{theorem}

The rest of the section is devoted to the proof of theorem
\ref{main.nf}. 

In order to put the system in normal form we will use
the method of Lie transform that we now recall.
Having fixed $r\geq 2$, consider a function $\chi$ of the form
\begin{align}
\label{chi}
\chi(M,\xi,\phi_c):=\sum_{{|\mu+\nu|=r}\atop {\omega\cdot(\nu-\mu)\not=0}}
\chi_{\mu\nu}(M)
\xi^\mu\bar\xi^\nu+\sum_{{|\mu|+|\nu|=r -1}\atop {|\omega\cdot(\nu-\mu)|<\Omega}}
\xi^{\mu}\bar\xi^{\nu}\left\langle E\Phi_{\mu\nu}(M); \phi_c
\right\rangle 
\end{align}
where $\chi_{\mu\nu}\in \ac(\R^n,\R)$, $\Phi_{\mu\nu}\in
\ac(\R^n,\W)$. 

\begin{remark}
\label{rem.ham}
If $\chi(M,\xi,\phi_c)$ is a smoothing Hamiltonian then, by lemma
\ref{l.6}, its Hamiltonian vector field generates a flow $\phi^t$ of
almost smoothing perturbations of the identity which are defined up to
time 1 in a sufficiently small neighborhood of the origin.
\end{remark}

Denote $\Tr:=\phi^1\equiv \phi^t\big|_{t=1}$. Such a transformation
will be called the {\it Lie transform} generated by $\chi$.

We define now the nonresonant projector $\Pinr$ acting on homogeneous
polynomials; it restricts the sum to nonresonant values of the
indexes. So let $F=F(M,\phi)$ be a homogeneous polynomial of degree
$r$ in $\phi$ continuous on $\V^j$. Consider first
$F_1(M,\phi):=F(M,\phi_d)+\di F(M,\phi_d)\phi_c$, where $\di$ is the
differential at fixed $M$,
\begin{align}
\label{chi.e}
F_1(M,\xi,\phi_c)=\sum_{|\mu+\nu|=r}F_{\mu\nu}(M)
\xi^\mu\bar\xi^\nu+\sum_{|\mu|+|\nu|=r -1}
\xi^{\mu}\bar\xi^{\nu}\left\langle E\Phi_{\mu\nu}^F(M); \phi_c
\right\rangle \ ,
\end{align}
in general $\Phi_{\mu\nu}^F\in \W^{-j}$, but we will see that in the
cases we will meet we have $\Phi_{\mu\nu}^F\in \ac(\R^n,\W)
$. We define
\begin{align}
\label{chi.F}
\Pinr F(M,\phi):=\sum_{{|\mu+\nu|=r}\atop {\omega\cdot(\nu-\mu)\not=0}}
F_{\mu\nu}(M)
\xi^\mu\bar\xi^\nu+\sum_{{|\mu|+|\nu|=r -1}\atop {|\omega\cdot(\nu-\mu)|<\Omega}}
\xi^{\mu}\bar\xi^{\nu}\left\langle E\Phi^F_{\mu\nu}(M); \phi_c
\right\rangle \ .
\end{align}

Finally, given $F(M,\xi,\phi_c)$, $F\in \ac(\R^n\oplus \V;\R)$, we
define the projector $\Pinr^r$ which by definition produces the
nonresonant part of the homogeneous Taylor polynomial of degree $r$ of
$F$.

In order to prove theorem \ref{main.nf} we proceed iteratively: we
assume the system to be in normal form at order $r-1$ and we normalize
it at order $r$. In order to perform the $r-th$ step we look for a
function $\chi_r(M,\phi)$, $\chi_r \in C_R(\Kcs,\R)$, $\Kcs=\R^n\oplus
\V$ such that the corresponding Lie transform $\Tr_r$ is the wanted
coordinate transformation.  Thus $\chi_r$ has to be chosen such that $
\Pinr^r\left( H\circ\Tr_r \right)=0$.  In order to write explicitly
such an equation, remark that, since the Lie transform generated by a
smoothing function is an almost smoothing perturbation of the
identity, after any number of coordinate transformations the
Hamiltonian has the form \eqref{H.new}-\eqref{H.new2}.

In order to compute the coefficients to be put equal to zero we work
in $\V^{k}$, with a sufficiently large $k$, so that also an almost
smooth map can be expanded in Taylor series up to order $r+1$.  Given
two functions $\chi(M,\phi)$ and $F(M,\phi)$ we denote by $
\left\{\chi;F\right\}^{st} $ the Poisson bracket of the two functions
           {\it computed at constant $M$, i.e. as if $M$ where
             independent of $\phi$}.  Similarly we will denote by
           $X_{\chi}^{st}(M,\phi)$ the Hamiltonian vector field of
           $\chi$ computed as if $M$ where independent of $\phi$.

We first study the simpler case in which $r\geq 3$.

\begin{lemma}
\label{lem.nf.1}
Assume that $\chi\in C_R(\Kcs,\R)$, $\Kcs:=\R^n\oplus
\V$, is a homogeneous polynomial of degree $r\geq3$, then one
has
\begin{equation}
\label{eq.nf.51}
\Pinr^r\left( H\circ\Tr_r \right)=\Pinr\left[\left\{H_{L0}+H_{L1};\chi
  \right\}^{st }+\frac{\partial D}{\partial M_j}(M) \langle
  A_j\phi;X^{st}_{\chi} \rangle\right]   +\Pinr^r H_N
\end{equation}
Furthermore one has $\Pinr^{r_1}(H\circ\Tr_r)=\Pinr^{r_1}(H)$
$\forall r_1<r$.
\end{lemma}
\proof 
First remark that, by Remark \ref{tei.r} and Remark \ref{rem.ham}, one
has ,
$$
\Tr_r(\phi)=\phi+X_{\chi}^{st}(\phi)+\Oc(|\phi|^{r})\ ,
$$ 
and therefore $H_N\circ\Tr_r=H_N+\Oc(|\phi|^{r+1})$, which shows that
$\Pinr^{r_1}\left( H_N\circ\Tr_r \right)=\Pinr^{r_1}\left( H_N \right)$
$\forall r_1\leq r$.
We come to $H_L\circ\Tr_r$. Denote $\phi'=\Tr_r(\phi)$ and 
$$
M'_j=M_j\circ \Tr_r= M_j+\langle
  A_j\phi;X^{st}_{\chi}(M,\phi) \rangle+\Oc(|\phi|^{r+1})\ ,
$$
then one has
\begin{align}
\label{hlcirctr}
H_L( & M',\phi')
= H_{L0}(\phi)+
\di H_{L0}(\phi)X_{\chi}^{st}(M,\phi) +\Oc(|\phi|^{r+1}) 
\\ 
\nonumber
&+
H_{L1}(M,\phi')+\Oc(|\phi|^{r+1})+ D(M)
\\
\nonumber
&+ \frac{\partial D}{\partial
  M_j}(M) \langle A_j\phi;X^{st}_{\chi}(M,\phi) \rangle
+\Oc(|\phi|^{r+1}) 
\\ 
\nonumber
&= H_{L0}(\phi)+ \di
H_{L0}(\phi)X_{\chi}^{st}(M,\phi)+ \di H_{L1}(M,\phi)X_\chi^{st}(\phi)
\\
\nonumber
&
+
D(M)+ \frac{\partial D}{\partial M_j}(M) \langle
A_j\phi;X^{st}_{\chi}(M,\phi) \rangle +\Oc(|\phi|^{r+1})
\end{align}
where the differentials are computed at constant $M$.  The application
of $\Pinr^r$ and the remark that $\di
H_{L1}X_{\chi}^{st}=\{H_{L1};\chi \}^{st}$ give the result.  
\qed

\begin{lemma}
\label{lem.nr.nf}
The function $\Pinr^r H_N$ is smoothing, and thus admits the
representation \eqref{chi.e} with $F_{\mu\nu}\in\ac(\R^n,\R)$ and
$\Phi_{\mu\nu}\in\ac(\R^n,\W)$.
\end{lemma}
\proof Consider $H_N(\phi_d+\phi_c)$, and remark that only
$H_N(\phi_d)$ and $\di H_N(\phi_d)\phi_c$ contribute to $\Pinr
H_N$. Now, $H_N(\phi_d)$ is clearly smoothing. To compute the other term
we use formula \eqref{H.new2}. Clearly the only term to be discussed
is the one coming from $H^3_P$. In order to compute it we use the
definition \eqref{e.ham.8.1} of $H^3_P$. Compute first the
differential with respect to $\phi$: 
\begin{equation}
\label{eq.nf.1}
\di H^3_P(\eta,\phi_d)\Phi=\di H_P (\eta+\phi_d)\Phi- \left[\di H_P
  (\eta)\Phi+\di ^2 H_P(\eta)(\phi_d,\Phi)  \right]\ ;
\end{equation} 
substituting $\Phi=\phi_c+\di S(\phi_d)\phi_c$ one gets the formula
for the term we have to compute. Here we denoted
$S:=S^1_1+S^0_2$. Then, by assumption (P2) the expression
\eqref{eq.nf.1} is smoothing.  \qed

\begin{lemma}
\label{homor>1}
For any $r\geq 3$ there exists $\chi_r\in C_R(\Kcs,\R)$,
$\Kcs=\R^n\oplus \V$, such that, denoting by $\Tr_r$ the
corresponding Lie transform, one has $\Pinr^r(H\circ\Tr_r)=0$.
\end{lemma}
\proof 
In order to construct $\chi_r$ we work with the implicit function
theorem, so it is useful to work in a fixed space $\W^k$ and with a
fixed regularity, say to work in $C^s$. For fixed $M$, we identify the
space of the functions of the form \eqref{chi} with the linear space
of its coefficients $(\chi_{\mu\nu};\Phi_{\mu\nu})$. When endowed by
the norm
\begin{equation}
\label{norchi}
\norma{\chi}_{k}:=\sup_{\mu\nu}
\left|{\chi_{\mu\nu}}\right|+\sup_{\mu\nu}
\norma{\Phi_{\mu\nu}}_{\W^{k}} \ ,
\end{equation}
it will be denoted by $\gen_{k}$.

We fix a suitable opens set $\U\subset \R^n$ and study the map
\begin{align}
\label{eq.r}
C^s(\U;\gen_{k})\ni \chi\mapsto \Le\chi:=\Pinr^r\left( H_L\circ\Tr_r
\right)
\\
\nonumber
=\Pinr\left[\left\{H_{L0}+H_{L1};\chi \right\}^{st
  }+\frac{\partial D}{\partial M_j}(M) \langle A_j\phi;X^{st}_{\chi}
  \rangle\right]
\in C^s(\U;\gen_{k-d_0}) \ .  
\end{align}
Remark that in \eqref{eq.r}, $M$ plays the role of a parameter, so we
work at fixed $M$ and consider 
$$
\Le:\gen_k\mapsto\gen_{k-d_0}\ .
$$ Such an operator is a relatively bounded perturbation of the linear
operator
\begin{equation}
\label{aab}
\gen_k\ni\chi\mapsto \Le_0\chi:=\Pinr \left\{
H_{L0};\chi\right\}^{st}\in\gen_{k-d_0}\ ,
\end{equation}
furthermore $\Le-\Le_0$ has a norm of order $M$. So we study
$\Le_0^{-1}$. To this end remark that one has
\begin{align}
\label{lochi}
\Le_0\chi=\sum_{\mu\nu}-\im \omega\cdot(\mu-\nu
)\chi_{\mu\nu}\xi^\mu\bar\xi^\nu 
\\
+
\sum_{\mu\nu}\left[\im \omega\cdot(\nu-\mu)\langle
  E\Phi_{\mu\nu};\phi_c \rangle-\langle EL_c\Phi_{\mu\nu};\phi_c
  \rangle \right] \xi^\mu\bar \xi^\nu\ ,
\end{align}
so that one has that
$\chi=\Le_0^{-1}\Pinr F$ is given by
\begin{equation}
\label{solhomosta4}
\chi_{\mu\nu}=\frac{F_{\mu\nu}}{-\im\omega\cdot
  (\mu-\nu)}\ ,\quad \text{for}\ \omega\cdot (\mu-\nu)\not=0
\end{equation}
\begin{equation}
\label{solhomosta3}
\Phi_{\mu\nu}=R_{L_c}(\im \omega\cdot (\nu-\mu))\Phi_{\mu\nu}^{F}\ .
\end{equation}
By (L1), the resolvent maps $\W^{k-d_0}$ into $\W^k$, thus it is
regularizing, therefore equations \eqref{solhomosta4},
\eqref{solhomosta3} show that the inverse of $\Le_0$ is smooth as a
map from $\gen_{k-d_0} $ to $\gen_{k} $. So $\Le$ can be inverted by
Neumann formula, giving $\chi$. The $C^s$ smoothness in $M$ follows
from the $C^s$ smoothness of of
$M\mapsto(F_{\mu\nu}(M),\Phi^F_{\mu\nu}(M))$. The thesis follows from
the arbitrariness of $k$ and $s$. \qed

We come now to the more complicated case $r=2$.  One has that
$X^{st}_\chi(\phi)=T(M)\phi$
with $T(M)$ a suitable linear smoothing operator.
Furthermore, remark that, by the proof of lemma \ref{l.6}, one has
\begin{equation}
\label{phi.0.e}
\phi'=\Tr_2(\phi)= \e^{T(M)}\phi+\Oc(|\phi|^{3})\Longrightarrow 
\Pinr^2(H_L\circ\Tr_2)= \Pinr^2(H_L\circ \e^{T(M)})\ .
\end{equation}
Thus
\begin{align}
\label{eq.hom.1}
\Pinr^2\left(H_L\circ\Tr_2\right)=
\Pinr^2\left(
\left\{\chi_2;H_{L0}\right\}^{st}
+H_{L1}(M,\phi)\right.
\\
\label{eq.hom.2}
+H_{L1}(M(\e^{T(M)}\phi),\e^{T(M)}\phi)-H_{L1}(M,\phi)
+D(M(\e^{T(M)}\phi))-D(M)
\\
\label{eq.hom.4}
\left.+
H_{L0}(\e^{T(M)}\phi)-\left[H_{L0}(\phi)+\left\{\chi_0;H_{L0}\right\}^{st}
  \right] 
\right)
\end{align}
which is a small perturbation of the first line. We will solve
$\Pinr^2\left(H_L\circ\Tr_2\right)=0$ using the implicit function
theorem, working perturbatively with respect to the first line.
First we need to estimate the other lines. To this end we need the
following lemma. 
\begin{lemma}
\label{lem.nor.0}
On the space of the functions $\chi\in\gen_{k}$ homogeneous of degree
2, the norm \eqref{norchi} is equivalent to the norm of
$X^{st}_{\chi}$ as a linear operator from $\Hc^{-k}$ to $\Hc^k$.
\end{lemma}
\proof One has
\begin{equation}
\label{eq.r.0}
X^{st}_\chi\left(
\begin{matrix}
\xi_k
\\
\phi_c
\end{matrix}
\right)= \left(
\begin{matrix}
-\im\sum_{|\mu|+|\nu|=2}\chi_{\mu\nu}\nu_k \frac{\xi^\mu\bar \xi^\nu}{\bar
  \xi_k} +\sum_{ |\mu|+|\nu|=1 }\nu_k \frac{\xi^\mu\bar \xi^\nu}{\bar
  \xi_k} \left\langle E\Phi_{\mu\nu},\phi_c  \right\rangle 
\\
\sum_{|\mu|+|\nu|=1}  \xi^\mu\bar\xi^\nu  \Phi_{\mu\nu}\ . 
 \end{matrix}
\right)
\end{equation}
so it is clear that the norm of such a linear operator is controlled
by the norm \eqref{norchi}. We have also to prove that the norm of
a single function $\Phi_{\mu\nu}$ (and the modulus of the coefficients
$\chi_{\mu\nu}$) is controlled by the norm of the linear operator. To
see this, remark that for example taking $\phi_c=0$ and $\xi_1=1$, $\xi_k=0$
for $k\not=1$, one gets
$T(\xi;\phi_c)=(\chi_{\mu^1\nu^1},\Phi_{\mu^2\nu^2}) $ with
$\mu^1=(1,0,0,0...)$ and so on. Thus the norm of each of the two
components is controlled by the operator norm of $X^{st}_{\chi}$. \qed

\begin{lemma}
\label{r=0}
There exists a $\chi_2\in C_R(\Kcs,\R)$, $\Kcs=\R^n\oplus \V$ of
the form \eqref{chi} with $r=2$, such that $\Pinr^2(H_L\circ\Tr_2)=0$
\end{lemma}
\proof Define 
\begin{equation}
\label{eq.o.hom}
G_1(M,\chi):=\Pinr\left(\eqref{eq.hom.2}+\eqref{eq.hom.4}
\right) 
\end{equation}
equation $\Pinr^2(H_L\circ \Tr_2)=0$ coincides with 
\begin{equation}
\label{eq.0.hom.1}
0=\Le_0\chi+\Pinr H_{L1}+G_1(M,\chi)\ ,\quad \Le_0\chi:=\left\{H_{L0};\chi
\right\}\ .
\end{equation}
Since by the same reasoning of the proof of lemma \ref{homor>1},
$\Le_0^{-1}$ is bounded as an operator from $\gen_{k} $ to
$\gen_{k+d_0}$, and the norm of $G_1(M,\chi)$ is smaller then
$C|M|\norma{\chi}$, one can apply the implicit function theorem to
\eqref{eq.0.hom.1}, getting the result. \qed

This concludes the proof of theorem \ref{main.nf}.

\section{Dispersive Estimates}\label{DisNLS} 

From now on we restrict our setting to the situation of NLS, but we
try to write clearly the assumptions we use, in order to make easy the
application to other models. Thus, from now on the scale $\Hc^k$ will
be that of the weighted Sobolev space $H^{k,l}$ (where the measure is
$\left\langle x\right\rangle^ld x$) and $\Hc^\infty=$Schwartz
space. When we write only one index we mean the standard Sobolev
spaces without weight. We will also use the Lebesgue spaces $L^p$ and
assume $d_0=2$.

In this section we will systematically use the notation $a\sleq b$ to
mean ``there exists a positive $C$, independent of all the relevant
quantities, s.t. $a\leq Cb$''.

Given functions $w^j(.)\in C^0([0,T],\R^n)$, consider the time
dependent linear operator 
\begin{equation}
\label{lin.time}
L(t):= P_cJ\left[A_0+V_0-w^j(t)A_j -\lambda^j(p_0)A_j\right]
P_c=L_c-w^j(t)J\A_j \ .
\end{equation}
and denote by $\U(t,s)$ the evolution operator of the equation
$\dot \phi=L(t)\phi$;
we assume that there exists $\epsilon>0$ such that, if
$|w^j(t)|<\epsilon$ then the following Strichartz  estimates
hold
\begin{itemize}
\item[(St.1)] \begin{align}
\label{flow}
\norma{\U(t,0)\phi_c}_{L^2_tL^6_x}\sleq \norma{\phi_c}_{L^2_x}\ ,
\\
\label{flow1}
\norma{\int_0^t\U(t,s)F(s)ds}_{L^2_tL^6_x}\sleq
\norma{F}_{L^2_tL^{6/5}_x}\ .
\end{align}
\item[(St.2)] There exists $a$ s.t., for any $k\geq0$ and any
  $\rho\in\pm\im(\Omega,\infty)$ the limit
\begin{equation}
\label{ris.l.e}
R^{\pm}_{L_c}(\rho):=\lim_{\epsilon\to0^+}(L_c-\rho\pm\epsilon)^{-1}\quad
\text{exists\ in}\ B(P_cH^{k,a},H^{k,-a})\ ;
\end{equation}
furthermore, for any $k, a\geq0$ one has 
\begin{equation}
\label{ris.l.e1}
[JV_0,JA_j]:H^{k,-a}\to H^{k,a}\ . 
\end{equation}

\item[(St.3)] For any $\Phi\in \W^\infty$, any
  complex valued function $h(.)\in L_t^2$, any $\rho \in \pm \im
  (\Omega,\infty)$ one has
\begin{align}
\label{flow2}
\norma{\langle x\rangle^{-a}\U(t,0)R^{\pm}_{L_c}(\rho)\Phi}_{L^2_x}\sleq
\frac{\norma {\langle x\rangle^a \Phi}_{L^2_x}}{\langle
  t\rangle^{3/2}}\ ,
\\
\label{flow3}
\norma{\int_0^t\U(t,s)h(s)R^{\pm}_{L_c}(\rho)\Phi ds}_{L^2_tL^{2,-a}_x}\sleq
\norma{h}_{L^2_t}\norma{\Phi}_{L^{2,a}_x}\ .
\end{align}
where $a$ is a sufficiently large constant.
\end{itemize}

Finally we need some smoothness of the vector field of $H_P$

\begin{itemize}
\item[(P3)] The map $X^2_P$ defined in \eqref{vf.31} is
  continuous as a map 
$$H^{k,l}\times (H^1\cap L^{6})\ni (\eta,\phi)\mapsto X^2_P\in
  L^{6/5}\ ,$$ with $k,l$ sufficiently large; furthermore, for
  $\phi\in H^1$, with $\norma{\phi}_{H^1}\leq \epsilon$, one has
\begin{equation}
\label{D2}
\norma{X^2_P(\eta;\phi)}_{L^{6/5}}\leq
\epsilon C  \norma{\phi}_{L^6} \ ,
\end{equation}
with $C=C(\epsilon,\norma{\eta}_{H^{k,l}})$. 

The map 
\begin{align}
H^{k,l}\ni \eta \mapsto
\di X_P(\eta)\in B(L^6,L^{6/5})
\end{align} 
is $C^1$ for large enough $k,l$.
\end{itemize}

The main result of this section is the following theorem
\begin{theorem}
\label{main.dis}
Consider the Hamiltonian \eqref{H.new} and assume it is in normal form
at order $2r_t$. Assume also that the Fermi Golden Rule \eqref{m.34} below
holds. Let $\phi(t)$ be a solution of the corresponding Hamilton
equations with an initial datum $\phi_0$ fulfilling
\begin{equation}
\label{eq.main.dis}
\norma{\phi_0}_{H^1}\leq \epsilon
\end{equation}
and $\epsilon$ small enough, then one has
\begin{align}
\label{eq.main.dis.1}
\norma{\phi_c(t)}_{L^2_tL^6_x}\sleq \epsilon
\\
\label{eq.main.dis.2}
\omega\cdot
\mu>\Omega\ ,\Longrightarrow\ \norma{\xi^\mu(.)}_{L^2_t}\sleq \epsilon
\ .
\end{align}
\end{theorem}

The rest of the section will be devoted to the proof of this
theorem. 

\subsection{Estimate of the continuous variable}\label{est.cont1}

Given a Hamiltonian of the form \eqref{H.new}, in
normal form at order $2r_t$ we study the solution of the corresponding
Hamilton equations. It will be denoted by $\phi(t)$.
\begin{remark}
\label{lem.F} 
Let $G$ be a map of the form $G(\phi)=\phi+S(\phi)$, with a smoothing
$S$. Consider the Hamiltonian $
H_P^3(\eta,G(\phi))$, with a fixed $\eta$. Then one has 
\begin{equation}
\label{jdfhi}
J\nabla(H_P^3\circ
G)=J\left[dG(\phi)\right]^*EX^2_P(\eta;G(\phi))\ ,
\end{equation}
where $ \di G(\phi)^* $ is the adjoint of the operator $\di G(\phi)$.
\end{remark}
\begin{remark}
\label{esti.M}
Using orbital stability (which
follows from (L2)) one has that, given an initial datum with $\norma
{\phi}_{H^{1}}\leq C\epsilon$, then
$$
|M_j(t)|\leq C\epsilon^2\ ,\quad \forall t\ .
$$
\end{remark}

\begin{remark}
\label{xp}
One has 
\begin{align}
\label{dd}
X_P^2(\eta;G(\phi_d&+\phi_c))-X_P^2(\eta;G(\phi_d))
\\
\nonumber
&=X^2_P(\eta+\phi_d+S(\phi_d);\phi_c+S(\phi_d+\phi_c)-S(\phi_d))
\\
\nonumber
&+\left[ \di X_P(\eta+\phi_d+S(\phi_d))-\di X_P(\eta)\right]
(\phi_c+S(\phi_d+\phi_c)-S(\phi_d))\ .   
\end{align}
\end{remark}

\begin{lemma}
\label{lem.nonlin}
Let $G$ be as in remark \ref{lem.F}, fix $\eta\in H^{k,l}$ with
sufficiently large $k,l$, and
consider a Hamiltonian function of the form
$H^3_P(\eta;G(\phi))+R^0_3$. Assume it is in normal form at order
$2r_t$. Let $X$ be its Hamiltonian vector field.
Assume that for some $T>0$ the functions $\phi_c(t)$ and $\xi(t)$ fulfill
the estimates
\begin{align}
\label{c.1}
\norma{\phi_c}_{L^2_t[0,T]L^6_x}\leq C_1\epsilon\ ,
\\
\label{c.2}
\norma{\xi^\mu}_{L^2_t[0,T]}\leq C_2 \epsilon\ ,\quad \forall \mu\in
\Kc:=\{ \mu\ :\ \omega\cdot \mu>\Omega\ ,\ |\mu|\leq 2r_t\}\ ,
\\
\label{c.311}
\norma {\phi_c(t)}_{H^1_x}\sleq \epsilon\ ,\quad |\xi(t)|\sleq
\epsilon\ ;
\end{align}
then there exists $C$ independent of $\epsilon, C_1,C_2$ such that one has
\begin{equation}
\label{c.3}
\norma{P_cX(\phi)}_{L^2_t[0,T]L^{6/5}_x}\leq \epsilon
C(C_2+\epsilon C_1)\ .
\end{equation}
\end{lemma}
\proof Write
$X(\phi_c+\phi_d)=X(\phi_d)+X(\phi_c+\phi_d)-X(\phi_d)$. From (P2),
\eqref{jdfhi}, \eqref{vf.31}, it
is clear that $X(\phi_d(\xi,\bar\xi))\in\ac(\C^K,\V)$.  Furthermore, write
$P_cX(\phi_d)=P_cX_{\leq r_t}(\phi_d)+P_cX_{> r_t}(\phi_d)$, where $P_cX_{\leq
  r_t}(\phi_d)$ is the Taylor expansion truncated at order $2r_t$,
which therefore (due to the fact that the system is in normal form)
contains only monomials of the form $\Phi_{\mu\nu}(M)\xi^\mu\bar
\xi^{\nu}$ with $\Phi_{\mu\nu}\in\ac(\R^n,\Wc)$ and
$|\omega\cdot(\mu-\nu)|>\Omega$.  This implies in particular
\begin{equation}
\label{mu.1}
\norma{P_cX_{\leq r_t}(\phi_d)}_{L^{6/5}_x}\sleq\sum_{|\mu|\leq 2r_t\ ,
  |\omega\cdot \mu|>\Omega} 
|\xi^\mu| \ .
\end{equation}
It follows from \eqref{c.2} that 
$$
\norma {P_cX_{\leq r_t}(\phi_d) } _{L^2_t[0,T]L^{6/5}_x}\leq CC_2\epsilon\ .
$$
Concerning $X_{> r_t}$, just remark that, by the formula for the
remainder of the Taylor expansion, one has
$$
\norma{P_cX_{> r_t}(\phi_d)}_{L^{6/5}_x}\sleq \left(
|\xi_1|^2+...+|\xi_K|^2 \right)^{(2r_t+1)/2}\ .
$$ Controlling the r.h.s. through $\norma{\xi^\mu}_{L^2_t}$, $\mu\in
\Kc$, one gets the thesis

We have now to estimate $P_cX(\phi_c+\phi_d)-P_cX(\phi_d)$. By remark
\ref{lem.F}, it is the sum of a smoothing term coming from $R_3^0$ and
of the quantity
\begin{align}
\label{x11}
[\di G^*(\phi_d+\phi_c)-\di
  G^*(\phi_d)]EX^2_P(\eta_{p_0-M};G(\phi_d+\phi_c))
\\
\label{x12}
+\di
G^*(\phi_d)E\left[X^2_P(\eta;G(\phi_d+\phi_c))-X^2_P(\eta;G(\phi_d))
  \right] \ .
\end{align}
Since $\di G^*(\phi_d+\phi_c)-\di
  G^*(\phi_d)=\di S^*(\phi_d+\phi_c)-\di
  S^*(\phi_d)$, (where we used the notations of lemma \ref{lem.F}) one has
$$
\norma{\di G^*(\phi_d+\phi_c)-\di
  G^*(\phi_d)}_{B(H^{-k_1,-l_1};H^{k_2,l_2})}\sleq
\norma{\phi_c}_{H^{-k_1,-l_1}}\sleq \norma{\phi_c}_{L^6}\ ,
$$
and therefore, 
\begin{align}
\label{x.13}
\norma{[\di G^*(\phi_d+\phi_c)-\di
    G^*(\phi_d)]EX^2_P(\eta_{p_0-M};G(\phi_d+\phi_c))}_{L^{6/5}}
\\ \nonumber \sleq
\norma{\phi_c}_{L^6}\norma{X^2_P(\eta;G(\phi_c+\phi_d))}_{L^{6/5}}
\sleq
\norma{\phi_c}_{L^{6}}\epsilon \norma{
  G(\phi_c+\phi_d) }_{L^6}
\\ \nonumber \sleq
 \epsilon
\norma{\phi_c}_{L^6}\norma{\phi}_{L^6}\sleq \epsilon
\norma{\phi_c}_{L^6}\norma{\phi}_{H^1} \sleq \epsilon^2\norma{\phi_c}_{L^6}\ .
\end{align}
In order to estimate \eqref{x12}, we exploit Remark \ref{xp} which
gives
\begin{align}
\label{x.15}
\norma{X^2_P(\eta;G(\phi_d+\phi_c))-X^2_P(\eta;G(\phi_d))}_{L^{6/5}}
\\
\nonumber
\sleq\epsilon
\norma{\phi_c+S(\phi_c+\phi_d)-S(\phi_d)}_{L^{6}}+\norma{G(\phi_d)}_{H^{k,l}}
\norma{\phi_c}_{L^6} \sleq \epsilon \norma{\phi_c}_{L^6}\ .
\end{align}
Adding the trivial estimate of $\di G^*(\phi_d)$ one gets the thesis.\qed

We are now ready to give the estimate of the continuous variable
$\phi_c$. 

\begin{lemma}
\label{lem.cont}
Let $\phi(t)$ be a solution of the considered system. Assume that the
initial datum $\phi$ fulfills $\norma \phi_{H^1}\leq \epsilon$ for
some $\epsilon$ small enough. Assume that, for some $T>0$, the a
priori estimates \eqref{c.1}, \eqref{c.2} and \eqref{c.311} are
fulfilled then $\phi_c$ fulfills the following estimate
\begin{equation}
\label{le.msu}
\norma{\phi_c(t)}_{L^2_t[0,T]L^6_x}\leq C\epsilon(C_2+\epsilon C_1)\ .
\end{equation}
\end{lemma}
\proof First, the equation for $\phi_c$ has the form
\begin{equation}
\label{fice}
\dot \phi_c=L(t)\phi_c+P_cJ[(V_M-V_0)\phi+ (S^1_1)_{lin}\phi]+P_c
X_N(\phi)\ .
\end{equation}
We also denoted by $(S^1_1)_{lin}$ a linear smoothing operator whose
norm tends to zero when $M\to0$. Remark that, since the system is in
normal form at order $r_t$, the term in square bracket is independent
of $\phi_d$. Thus \eqref{fice} is equivalent to
\begin{equation}
\label{fice1}
\dot \phi_c=L(t)\phi_c+P_cJ[(V_M-V_0)\phi_c+ (S^1_1)_{lin}\phi_c]+P_c
X_N(\phi)\ .
\end{equation}
We use Duhamel principle to write its solution in the form
$\phi_c(t)=I_1+I_2+I_3+I_4$, 
where
\begin{align}
\label{if.1}
I_1:=\U(t,0)\phi_c(0)\ ,
\quad I_2:= \int_0^t \U(t,s)P_cJ(V_M-V_0)\phi_c(s) \di s\ ,
\\
\label{if.3}
I_3:=\int_0^t \U(t,s)P_c(S_1^1)_{lin}\phi_c(s) \di s\ ,
\quad
I_4:=\int_0^t \U(t,s)P_c X_N(\phi(s))\di s\ .
\end{align}
The estimates of $I_1$ and of $I_4$ are an immediate consequence of
(St.1) and lemma \ref{lem.nonlin}, which give
\begin{equation*}
\norma{I_1}_{L^2_t[0,T]L^6_x}\sleq \epsilon\ ,\quad
\norma{I_4}_{L^2_t[0,T]L^6_x}\sleq \epsilon(C_2+\epsilon C_1)\ .\quad 
\end{equation*}
Concerning $I_2$ we have, by the second of (P3)
$$
\norma{ (V_M-V_0)\phi_c(s)}_{L^{6/5}_x}\sleq |M|\norma{\phi_c}_{L^{6}_x}\sleq
\epsilon^2\norma{\phi_c}_{L^6_x}\ ,
$$
from which $\norma{I_2}_{L^2_t[0,T]L^6_x}\sleq
\epsilon^3C_1$. Similarly, $I_3$ is estimated using 
$$
\norma{P_c(S^1_1(M))_{lin}\phi_c}_{L^{6/5}_x}\sleq
\norma{(S^1_1(M))_{lin}\phi_c}_{H^{k,l}}\sleq |M|
\norma{\phi_c}_{H^{-k,-l}}\sleq \epsilon^2\norma{\phi_c}_{L^6}\ .  
$$
from which $\norma{I_3}_{L^2_t[0,T]L^6_x}\sleq \epsilon^3C_1$, and the
thesis. \qed

\subsection{A further step of normalization}\label{nor} 

Consider again the Hamiltonian in normal form at order $2r_t+1$, we
rewrite it in a form suitable for the forthcoming developments. First
write
\begin{equation}
\label{hnd}
H_{re}(\phi_d,\phi_c):=H_N(\phi_d+\phi_c)-H_N(\phi_d)-\di
H_N(\phi_d)\phi_c\ ,
\end{equation}
and
\begin{align}
\label{d0}
H_{Nd}:=H_N(\phi_d(\xi,\bar\xi))-Z_0(\xi,\bar \xi)\ ,
\quad
H_{Nc}:=\di
H_N(\phi_d)\phi_c- Z_1(\xi,\bar\xi,\phi_c)\ , 
\end{align}
where $Z_0$ is the Taylor expansion of $H_N(\phi_d(\xi,\bar\xi))$
truncated at order $2r_t+1$, and we defined
\begin{align}
\label{m.2} Z_1(\xi,\bar\xi,\phi_c):= \langle EG,\phi_c\rangle +\langle
E\bar G,\phi_c\rangle \ , 
\\ 
\label{m.2.21}
G:=\sum_{\nu\in\Kc}\bar \xi^\nu
\Phi_{\nu}\ ,\quad \bar G=\overline{\sum_{\nu\in\Kc }\bar
  \xi^\nu \Phi_{\nu}},
\end{align}
with $\Phi_{\nu}\in \ac(\R^n,\W)$.
Denote also
$\resto:=H_{Nd}+H_{Nc}+H_{re}$,
then the Hamilton equations of the system can be written in the form
\begin{align}
\label{m.4} 
\dot \phi_c &=L_c\phi_c+J\nabla_{\phi_c}H_{L1}+G+\bar G+
w^j(M,\xi,\phi_c) J\A_j\phi_c+J \nabla_{\phi_c}\resto
(M,\xi,\phi_c)\ , 
\\
\label{m.5} \dot \xi_k&= -\im  \omega_k \xi_k- \im \frac{\partial
H_{L1}}{\partial
  \bar \xi_k}- \im \frac{\partial
Z_0}{\partial
  \bar \xi_k}-\im \left\langle E\frac{\partial G}{\partial
  \bar \xi_k},\phi_c\right\rangle 
-\im \frac{\partial
  \resto}{\partial \bar \xi_k} \ ,
\\
&w^j(M,\xi,\phi_c):=\frac{\partial H_N}{\partial M_j}+\frac{\partial
  H_{L1}}{\partial M_j}+\frac{\partial D}{\partial M_j} \ , 
\end{align}
and the gradient $\nabla_{\phi_c}$ is computed at constant $M$. 
We look now for functions
$Y_\nu=Y_\nu(M)$ such that the new variable $g$ defined by
\begin{equation}
\label{m.6}
g:=\phi_c+Y+\bar Y\ ,\quad Y=\sum_{\nu}Y_\nu\bar \xi^\nu
\end{equation}
is decoupled up to higher order terms from the discrete variables.
Substitution into equation \eqref{m.4} yields
\begin{align}
\label{m.7.a} \dot g=L_cg+\sum_\nu\left(\Phi_\nu+\im \nu\cdot \omega
Y_\nu-L_cY_\nu\right)\bar \xi^\nu
\\
+\sum_\nu\left(\bar \Phi_\nu-\im \nu\cdot \omega
\bar Y_\nu-L_c\bar Y_\nu\right) \xi^\nu
+{\rm h.o.t.}\ 
\end{align}
where the h.o.t. will be explicitly computed below. 
In order to kill the main terms define 
\begin{equation}
\label{m.11}  Y_{\nu} = R^{+}_{L_c}(\im \omega\cdot \nu) \Phi_{ \nu }
\quad \text{and} \quad \bar Y_{\nu} =\overline{R^{+}_{L_c}(\im
  \omega\cdot \nu)  
\Phi_{ \nu  }} =R^{-}_{L_c}(-\im \omega\cdot \nu) \bar\Phi_{ \nu }\ .
\end{equation}

We substitute \eqref{m.6} into \eqref{m.5}. Then, using \eqref{m.11},
we get
\begin{align}
\label{m.13} 
\dot \xi_k &=  -\im   \omega_k \xi_k- \im
\frac{\partial Z_0}{\partial \overline{\xi}_k}- \im \frac{\partial
  H_{L1}}{\partial \bar \xi_k}+\G_{0,k}(\xi)-\im \left\langle
E\frac{\partial G}{\partial\bar \xi_k};g\right\rangle -\im
\frac{\partial \resto}{\partial\bar \xi_k}\ , 
\\
\label{m.13ter}& \G_{0,k}(\xi):=\im     \sum_{   \nu  \in \Kc,  \nu' \in
\Kc }\left(  \frac{\xi^{ \nu' }\bar \xi^{ \nu}} { \bar\xi_k} \nu_kc_{\nu\nu'
}+
 \frac{\bar\xi^{ \nu' }\bar \xi^{ \nu}} { \bar\xi_k} \nu_kb_{\nu\nu'
}\right)\ ,
\\ 
&c_{\nu\nu'}:=\langle E\Phi_{\nu },\bar Y_{\nu'} \rangle  \ ,
b_{\nu\nu'}:=\langle E\Phi_{\nu },Y_{\nu'} \rangle  \ . 
\end{align}

Following the standard ideas of normal form theory, we look for a
change of variables of the form $
\eta _j =\xi _j + \Delta_j(\xi)$
which moves to higher order the nonresonant terms. The choice
\begin{align}
\label{m.17bis}
\Delta_j(\xi):= \sum_{\substack{\mu \in \Kc, \, \nu\in \Kc \\
\omega\cdot(\mu-\nu )\not=0}}\left[
 \frac{1}{\im \omega \cdot (\mu - \nu) }
\frac{\xi ^{ \mu }\bar {\xi }^ {\nu }}{\bar \xi_j} \nu _j c_{\nu\mu
}
+ \frac{1}{-\im \omega \cdot (\mu + \nu) }
\frac{\bar\xi ^{ \mu }\bar {\xi }^ {\nu }}{\bar \xi_j} \nu _j b_{\nu\mu
}
\right]
\end{align}
transforms \eqref{m.13} into the system $\dot \eta_k=\Xi_k(\eta,\bar
\eta)+\E_k(t)$ where
\begin{align}
\label{m.17.terr}
\Xi_k(\eta,\bar \eta):= -\im \omega_k
\eta_k- \im \frac{\partial Z_0}{\partial
  \bar \eta_k}- \im \frac{\partial
H_{L1}}{\partial
  \bar \xi_k}+\Nc_k(\eta)
\\
\label{m.17ter}
\Nc_k(\eta):=\im  \sum_{\substack{\mu \in \Kc, \, \nu\in \Kc \\
\omega\cdot(\mu-\nu ) =0}} \frac{\eta^{\mu}\bar\eta^{\nu}}
{\bar\eta_k} \nu _kc_{\nu\mu }\ ,
\end{align}
and $\E_j(t)$ is a remainder term whose expression is explicitly
given by
\begin{align}
\label{ej}
&\E_j := X_j^{L1}(\xi)-X_j^{L1}(\eta)+\G_{0j}(\xi)-\G_{0j}(\eta)+X_j^N
\\
\nonumber
&+\sum_k \left(\frac{\partial \Delta_j}{\partial\xi_k}(\xi)\left[
  X^{L1}_k+\G_{0k}(\xi) +X^N_k(\xi)  \right]
+\frac{\partial \Delta_j}{\partial\bar \xi_k}(\xi)\left[
  \bar X^{L1}_k+\bar \G_{0k}(\xi) +\bar X^N_k(\xi)  \right]\right)
\\
\nonumber
&-\sum_k\left( \frac{\partial \Delta_j}{\partial\xi_k}(\xi)\im \omega_k\xi_k 
-\frac{\partial \Delta_j}{\partial\xi_k}(\eta)\im \omega_k\eta_k 
+ \frac{\partial \Delta_j}{\partial\bar \xi_k}(\xi)\im
\omega_k\bar \xi_k 
- \frac{\partial \Delta_j}{\partial\bar  \xi_k}(\eta)\im
\omega_k\bar \eta_k\right) 
\ ,
\end{align}
and we denoted
\begin{align*}
X_{k}^{L1}:= -\im \frac{\partial Z_0}{\partial \bar\xi_k}-\im
\frac{\partial H_{L1}}{\partial \bar\xi_k}\ ,
\quad
X_k^N:= \im \frac{\partial \resto}{\partial \bar\xi_k}-\im
\left\langle E\frac{\partial G}{\partial \bar \xi_k};g\right\rangle   \ .
\end{align*}

The key point is that the considered system for $\eta$ is no more
conservative. To see this we compute the Lie derivative of
$H_{0L\xi}\equiv \sum_{k}\omega_k\left|\eta_k\right|^2$ with respect
to $\Xi$. 

We partition $\Kc$ into ``resonant sets''. Define 
$$
\Lambda:=\left\{\lambda\in\R\ :\ \lambda= \omega\cdot \mu\ ,\quad
\mu\in \Kc\right\}
$$ and, for $\lambda\in\Lambda$, define
\begin{align}
\Kc_\lambda:=\left\{ \mu\in
\Kc\ :\ \omega\cdot \mu=\lambda\right\} \ 
\\
\label{m.32.e}
 F_\lambda(\eta):= \sum_{\mu\in\Kc_\lambda}
\Phi_\mu\bar \eta^\mu\in\ac(\R^n,\W)\ .
\end{align}

\begin{remark}
\label{l.dis} One has
\begin{equation}
\label{m.33} \lie_{\Xi}H_{0L\xi}=-\Im\left( \sum_{\lambda\in\Lambda}\lambda
\left\langle E
\overline{F_\lambda(\eta)};R^+_{L_c}(\im \lambda){F_\lambda(\eta)}
\right\rangle\right) \ .
\end{equation}
Furthermore, using formally the formula $(x-\im 0)^{-1}=PV(1/x)+\im
\pi\delta(x)$ in order to compute (formally) $R^+_{L_c}(\im \lambda)$,
one realizes that if there are no convergence problems, one has
\begin{equation}
\label{fgr1}
\Im\langle E\bar \Phi;R^+_{L_c}(\im \lambda)\Phi\rangle >0\ ,\quad
\forall0\not=\Phi\in P_cH^{k,a} \ 
\end{equation}
with $k,a$ sufficiently large. 
In typical cases (e.g. in NLS), \eqref{fgr1} is obtained by using the
wave operator in order to conjugate $L_c$ and $JA_0$, and exploiting the result
by \cite{Yaj,Cuc01} according to which the wave operator leaves
invariant the $L^p$ spaces. 
\end{remark}

We are ready to state the Fermi Golden Rule, which essentially states
that the quantity \eqref{m.33} is nondegenerate; to this end denote
$$
b_\lambda(\eta):=\Im\left(\lambda
\left\langle E
\overline{F_\lambda(\eta)};R^+_{L_c}(\im\lambda){F_\lambda(\eta)}
\right\rangle\right)
$$
\begin{itemize}
\item[(FGR)] there exists a positive constant $C$ and a sufficiently
small $\delta _0>0$ such that for all $|\eta |<\delta _0$
\begin{equation}
\label{m.34}
\sum_{\lambda\in\Lambda}\lambda b_\lambda(\eta) 
\geq C {\sum_{\mu\in \Kc}|\eta^\mu|^2} \ .
\end{equation}
\end{itemize}
\begin{remark}
\label{wei}
This version of the FGR is essentially identical to that used in
\cite{GW08}. It is easy to see that in the nonresonant case
$\#\Kc_\lambda=1$ $\forall \lambda\in\Lambda$, \eqref{m.34} is
equivalent to the assumption that a finite number of coefficients is
different from zero (see \cite{BC09} condition (H7'')).
\end{remark}

\subsection{Estimate of the variables $g$, $\xi$, $\eta$.}

In order to estimate the variable $g$ we need the following
lemma
\begin{lemma}
\label{lem.commu}
For any $\Phi\in P_cH^{k,a}$ and any $\rho\in\sigma_c(L_c)$, there
exists a $\Psi\in P_cH^{k,a}$, linearly dependent on $\Phi$, such that
one has
\begin{equation}
\label{commu1}
\left[R^{\pm}_{L_c}(\rho),J\A_j \right] \Phi= R^{\pm}_{L_c}(\rho)\Psi\ .
\end{equation}
\end{lemma}
\proof To start with take $\rho\not \in \sigma(L_c)$. A simple
computation shows that (omitting $\rho$) the l.h.s. of \eqref{commu1}
is given by $R_{L_c}P_c[JV_0;JA_j] R_{L_c}\Phi$, from which, using (St.2)
and taking the limit $\rho\to\sigma_c$, the thesis follows.\qed

\begin{lemma}
\label{lem.cont.g}
Under the same assumptions of lemma \ref{lem.cont}, $g$ fulfills the
estimate
\begin{equation}
\label{esti.g}
\norma{g}_{L^2_t[0,T]L^{2,-a}_x}\leq C_0 \epsilon+ C\epsilon^2\ , 
\end{equation}
where $a$ is a sufficiently large constant and $C_0$ depends
only on the constant of the inequality \eqref{flow2}. 
\end{lemma}
\proof Remarking that (where $\nabla_{\phi_c}$ is computed at constant
$M$)
\begin{equation}
\label{remaf}
J\nabla_{\phi_c}H_{re}=P_c\left[X_N(\phi_c+\phi_d)-X_N(\phi_d)\right]\ ,\quad 
J\nabla_{\phi_c}H_{Nc}=P_cX(\phi_d)-(G+\bar G)\ ,
\end{equation}
and denoting $\resto_\xi:=H_{L1}+H_N$, the equation for $g$ takes
the form
\begin{align}
\label{eq.g.1}
\dot g&=L(t)g-w^jJ\A_j(Y+\bar Y) 
-\im \frac{\partial
  \bar Y}{\partial \xi_k}\frac{\partial \resto_\xi}{\partial
  \bar\xi_k}+ \im \frac{\partial Y}{\partial
  \bar\xi_k}\frac{\partial \resto_\xi}{\partial\xi_k}
\\
\label{eq.g.3}
&+J\nabla_{\phi_c} H_{L1}(\phi_c)+P_c\left[X_N(\phi_c+\phi_d)-X_N(\phi_d) \right]
\\
\label{eq.g.4}
&+P_c\left[ X(\phi_d)-(G+\bar G)\right] \ .
\end{align}
We apply Duhamel formula and estimate the different terms
arising. First we consider 
\begin{equation}
\label{eq.g.5}
\int_0^t\U(t,s)w^jJ\A_jY(s)\di s\ . 
\end{equation}
Using lemma \ref{lem.commu} and formula \eqref{m.11} for $Y$, it can be
rewritten as the sum of finitely many terms of the form
\begin{align}
\label{eq.g.51}
\int_0^t\U(t,s)w^jP_c
\left[ R^{\pm}_{L_c} J\A_j \Phi_\mu +R^{\pm}_{L_c}\Psi_\mu\right]\bar
\xi^\mu(s) \di s \ ,      
\end{align}
with suitable $\Psi_\mu\in\W^\infty$. Then, exploiting
(St.3\null), one has
\begin{equation}
\label{eq.5.52}
\norma{\eqref{eq.g.51}}_{L^2_t[0,T]L^{2,-\nu}_x}\sleq
\norma{\xi^{\nu}}_{L^2_t}\norma{\Phi_\nu}_{L^{2,\nu}_x}\sleq
C_2\epsilon^3  \ .
\end{equation}
The estimate of the last term of \eqref{eq.g.1} is identical to the
same estimate of \cite{BC09}, see lemma 7.9, so it is omitted. The
terms coming from \eqref{eq.g.3} have already been estimated in the
proof of lemma \ref{lem.cont} (see the estimates of $I_1,I_2,I_3$).

We come to \eqref{eq.g.4}. To this end remark that one has
\begin{equation}
\label{6.0}
X(\phi_d(\xi))=\sum_{\mu,\nu\ :\ \omega\cdot(\nu-\mu)>\Omega }
(X_{\mu\nu}\xi^\mu\bar \xi^\nu+c.c.)
\end{equation}
(with c.c. denoting the complex conjugated term), while the term
subtracted in \eqref{eq.g.4} coincides with
\begin{equation}
\label{6.1}
\sum_{\nu\in \Kc}
(X_{0\nu}\bar \xi^\nu+c.c.)\ .
\end{equation}
It follows that, if a term is present in \eqref{6.0} but not in
\eqref{6.1}, then it is of the form $X_{\mu\nu}\xi^\mu\bar \xi^{\nu+\nu'}$
with $\nu\in \Kc$. It follows that for such a term
\begin{equation}
\label{eq.g.7}
\norma{X_{\mu({\nu+\nu'})}\xi^\mu\bar \xi^{\nu+\nu'}}_{L^2_t[0,T]L^{2,-a}_x}\sleq
\norma{X_{\mu(\nu+\nu')}}_{L^{2,-a}_x}
\norma{\xi^\nu}_{L^2_t[0,T]}|\xi^\mu|\sleq C_2 \epsilon^2\ .
\end{equation}
Since the sum is finite the thesis follows. \qed

\begin{lemma}
\label{l.es.11}
Assume \eqref{c.1} and \eqref{c.2}, then, provided $\epsilon$ is
small enough, the following estimate holds
\begin{equation}
\label{e.es.117}
\sum_j\norma{\eta_j\E_j}_{L^1_t[0,T]}\leq CC_2\epsilon^2
\end{equation}
\end{lemma}
The proof of this lemma is almost identical to the proof of Lemma 7.11
of \cite{BC09}. Indeed the only difference is due to the presence of
$H_{L1}$, but the corresponding terms can be estimated by the same
methods used in \cite{BC09}. For this reason we omit the proof.

\begin{theorem}
\label{maintech}
Assume \eqref{c.1} and \eqref{c.2} then, provided $\epsilon$ is small
enough, one has
\begin{align}
\label{c.T.1}
\norma{\phi_c(t)}_{L^2_t[0,T]L^6_x}\leq C(C_2) \epsilon
\\
\label{c.T.2}
\omega\cdot\mu>\Omega\ \Longrightarrow\ \norma{\xi^\mu(t)}_{L^2_t[0,T]}\leq
C \sqrt{C_2} \epsilon
\end{align}
\end{theorem}
The proof (by standard bootstrap argument) is identical to the proof
of Theorem 7.12 of \cite{BC09} and therefore is omitted.

Then also Theorem \ref{main.dis} immediately follows.

\section{NLS}\label{NLS} 

Consider the scale of {real} Hilbert spaces $H^{k,l}(\R^n,\C)$. We
introduce the scalar product in $H^0$ and the symplectic form as
follows:
\begin{align}
\label{1.1}
\langle \psi_1;\psi_2\rangle:=2 \Re \left(\int_{\R^n}\psi_1(x)\bar
\psi_2(x)dx \right)\ ,\quad 
\omega(\psi_1;\psi_2 ):=\langle \im \psi_1;\psi_2\rangle\ ,
\end{align}
(remark that on a real vector space the multiplication by i is not a
scalar but a linear operator). The
Hamilton equations are given by $
\dot \psi=-\im \nabla _{\bar \psi}H$,
where $\nabla_{\bar\psi}$ is the gradient with respect of the $L^2$
scalar product.

The Hamiltonian of the NLS is given by
\begin{equation}
\label{1.6}
H:=\Phzero+H_P\ ,\qquad \Phzero(\psi):= \int_{\R^3}\bar\psi\,
(-\Delta\psi)d^3x\ ,\quad H_P(\psi):=-\int_{\R^3}\beta(|\psi|^2) d^3x \ ;
\end{equation}
in particular one has $A_0:=-\Delta$. The corresponding Hamilton
equations are \eqref{NLS1}.

There are 4 symmetries: Gauge and translations.  The operators
generating the symmetries are $A_j=-\im\partial_j $, $j=1,2,3$ and
$A_4=\uno$, so that one has
\begin{equation}
\label{sym.NLS.1}
\Ph_j(\psi)=-\int_{\R^n}\bar \psi \im\partial_j\psi d^nx\ ,\quad
\Ph_4(\psi)= \int_{\R^n}|\psi|^2 d^nx\ .
\end{equation}

The construction of the ground state (and the subsequent study)
exploits the boost transformation which, given a ground state at rest,
puts it in uniform motion. 
\begin{definition}
\label{boost}
Given $v\in\R^3$, the unitary transformation 
\begin{equation}
\label{boost.1}
U(v):L^2\ni \psi \mapsto U(v)\psi:=e^{-\frac{ i v\cdot x}{2}}\psi\in
L^2\ ,
\end{equation}
is called the {\it boost} transformation with velocity $v$.
\end{definition}
\begin{remark}
\label{boost.2}
The boosts form a unitary group parametrized by the
velocities. Furthermore, for any fixed $v$ the corresponding boost is
a canonical (symplectic) transformation. The boosts have also the
remarkable property of conserving the $L^p$ norms.
\end{remark}

Having fixed $\E >0$ and putting $\lambda^4=-\E $ and $\lambda^j=0$,
for $j=1,2,3$, the equation \eqref{e.ham.4} for the ground state,
denoted by $b_\E $, takes the form \eqref{ground.NLS.0} which has
already been discussed. We will
denote
$$
p_4(\E ):=\Ph_4(b_\E )\equiv\int_{\R^n}b_\E ^2d^nx\ .
$$

A direct computation shows that $\eta_p:=U(v)b_\E $ is a ground state
with parameters
\begin{align}
\label{grou.21}
\Ph_4(\eta_p)=p_4(\E )\ ,\quad \Ph_j(\eta_p)=\frac{v_j}{2}p_4(\E ),
\ j=1,2,3
\\
\lambda^j=v_j\ ,\ j=1,2,3\ ,\quad
\lambda^4=-\left(\E +\frac{|v|^2}{4}\right) \ .
\end{align}

%
%
%

In order to explicitly perform the computations and to verify all the
assumptions it is useful to exploit the existence of the boosts.  So
fix $p_0$ and consider $\eta_{p_0}$ and the decomposition of $H^k$
into $T_{\eta_{p_0}}\Tc$ and its symplectic orthogonal. Since $U(v)$
maps $b_\E $ into $\eta_{p_0}$, is linear and symplectic, it also maps
$T_{b_\E }\Tc$ to $T_{\eta_{p_0}}\Tc$ and $T_{b_\E }^\omega\Tc$ to
$T_{\eta_{p_0}}^\omega\Tc$. Furthermore it is unitary (and it also
conserves all the $L^p$ norms), and therefore it {\it is particularly
  convenient to represent $\V^k$ as $\V^k=U(v)\V^k_\E $, with $\V^k_\E :=
  T_{b_\E }^\omega\Tc\cap \Hc^k$. {\bf This is what we are now going to do}. We
  will also denote by $\Pi_\E \equiv \Pi_{b_\E }$ the projector on
  such a space}

\begin{remark}
\label{r.lin.nls}
In such a representation one has that $H_{L0}$ is represented by the
restriction to $\V_\E ^k$ of
$H_{L\E }:=\Phzero+d^2H_P(b_\E )(\psi,\psi)+\E \Ph_4$. Correspondingly the linear
operator $J\nabla H_{L0}$ is equivalent (through $U$) to the
restriction to $\V_\E ^k$ of the vector
field of $H_{LE}$, which in turn is the Hamiltonian vector field of
$H_{LE}\circ\Pi_\E $.   
\end{remark}

\begin{remark}
\label{r.lin.2}
Explicitly, one has
\begin{equation}
\label{e.lin.101}
H_{LE}(\psi)= \frac{1}{2}\left\langle -\Delta\psi,\psi
\right\rangle+\E \frac{1}{2} \langle \psi,\psi\rangle+
d^2H_P(b_\E )(\psi,\psi)\ , 
 \end{equation}
or, denoting
\begin{equation}
\label{1.8}
\psi=\frac{\psi_-+\im \psi_+}{\sqrt2}\ ,\quad \psi_{\pm}\in
H^{k,l}(\R^n,\R)\ ,
\end{equation}
\begin{equation}
\label{e.lin.11.2}
H_{LE}(\psi_+,\psi_-)= \frac{1}{2} \left\langle A_+\psi_+;\psi_+
\right\rangle + \frac{1}{2} \left\langle A_-\psi_-;\psi_-
\right\rangle\ ,
\end{equation}
where 
\begin{equation}
\label{e.lin.12.1}
A_-:=-\Delta+\E -\beta'(b_\E ^2)-2\beta''(b_\E ^2)b^2_\E \ ,\quad
A_+:=-\Delta+\E -\beta'(b_\E ^2)\ ,
\end{equation}
one has
\begin{equation}
\label{lzer0}
L_0\left[
\begin{matrix}
\psi_+
\\
\psi_-
\end{matrix}
\right]=
\left[
\begin{matrix}
-A_-\psi_-
\\
A_+\psi_+
\end{matrix}
\right] 
\end{equation}
\end{remark}

We pass to the verification of the assumptions. (S1-S4) are
trivial. The same is true for (P1-P2), (B1-B2). (L1,L2,P3) are well known
in this context, while (L3,L4) were assumed explicitly in
sect. \ref{StateNLS}. (St.2) is by now standard.  
We come to the other assumptions.

\begin{lemma}
\label{perelman.1}
Assumption (B3) holds.
\end{lemma}
\proof It is clearly enough to verify the assumption at $\eta_p=b_{\E}$. First
remark that, at $b_\E $, we have 
\begin{align}
\label{bE}
\frac{\partial \eta_p}{\partial p_j}=-\frac{\im
}{p_4}x^jb_\E \equiv\left(-\frac{1
}{p_4}x^jb_\E,0 \right)\ ,\ j=1,2,3\ ,
\\
\frac{\partial
  \eta_p}{\partial p_4} =-\E '\frac{\partial b_\E }{\partial \E }\equiv
\left(0, -\E '\frac{\partial b_\E }{\partial \E }\right)\ ,
\end{align}
so that one gets
\begin{align*}
\omega\left(\dep{\eta_p}4;\dep{\eta_p}j\right)=-\frac{\E'}{p_4}
\left\langle \left(\frac{\partial b_{\E}}{\partial E},0\right);
\left(x^jb_{\E},0\right)\right\rangle 
\\
= \int_{\R^3}x^jb_{\E}\frac{\partial b_{\E}}{\partial \E} d x\ ,
\end{align*}
but $b_{\E}\frac{\partial b_{\E}}{\partial \E}$ is spherically
symmetric, while $x^j$ is skew symmetric, and thus the integral
vanishes. \qed

The following
Lemma is a minor variant of a Lemma proved by Beceanu and Perelman
\begin{lemma}
\label{perelman}
Assumptions (St.1,St.3\null) hold.
\end{lemma}
In appendix \ref{B} we report its proof following \cite{appPer}.

\begin{corollary}
\label{cor.NLS}
Under the assumptions of section \ref{StateNLS}, Theorem \ref{main.dis}
holds for the NLS.
\end{corollary}

\begin{remark}
\label{pere.2}
In the case of NLS one can easily show that (St.1) and (St.3\null) hold
also if the spaces $L^p_x$ are substituted by the Sobolev spaces
$W^{1,p}_x$. As a consequence also the conclusion
\eqref{eq.main.dis.1} holds with $W^{1,6}_x$ in place of
$L^6_x$. 
\end{remark}

In the case of NLS, the flow of $L_0$ is well known to satisfy
Strichartz estimates of the form 
\begin{align}
\label{str.NLS.1}
\norma{\e^{tL_0}P_c\phi}_{L^q_tW^{1,r}_x}\sleq \norma{\phi}_{H^1}
\\
\label{str.NLS.3}
\norma{\int_0^t\e^{(t-s)L_0}P_cF(s)\di s}_{{L^{q}_tW^{1,r}_x}  }\sleq
\norma{F}_{L^{\tilde q'}_tW^{1,\tilde r'}_x}
\end{align}
for all admissible pair $(q,r)$, $(\tilde q,\tilde r)$, namely pairs
fulfilling
$$
2/q + 3/r = 3/2 ,\quad 6 \geq r \geq 2 , q\geq 2 \ .
$$
As a consequence one can prove the same estimates also for the flow
$\U(t,s)$. Using such estimates one gets the following 

\begin{theorem}
\label{NLS.scat}
Let $\phi(t)$ be a solution of the reduced equations corresponding to
NLS. Let $\phi_0$ be the initial datum, and assume
$\norma{\phi_0}_{H^1}=\epsilon$ is small enough. Then there exists
$\phi_\infty$ such that
\begin{equation}
\label{scat.NLS.finale}
\lim_{t\to+\infty}\norma{\phi(t)-e^{tL_0}\phi_\infty }_{H^1}=0
\end{equation}
\end{theorem}
\proof The proof is standard (see e.g. \cite{BC09}, Lemma 7.8) and
thus it is omitted.  Theorem \ref{main.NLS} is just a reformulation of
the above theorem in terms of the original system.

\appendix

\section{Proof of theorem \ref{t.red.1}}\label{red.proof}

We recall the idea on which the proof is based in the standard
case. The main point is the construction of a suitable coordinate
frame in which the actions of the symmetries becomes trivial. 

To start with consider the map
\begin{equation}
\label{m.cco.con}
I\times\R^n\times\V^k\ni(p,q,\phi)\mapsto
e^{q^jJA_j}(\eta_p+\pip\phi)\in \Hc^k\ .
\end{equation}

\begin{lemma}
\label{l.con.coo}
There exists a mapping $\varphi(u)\equiv (p(u),q(u))$ with the
following properties
\begin{itemize}
\item[1)] $\forall k$ there exists an open neighborhood
  $\U^{-k}\subset \Hc^{-k}$ of $\eta_{p_0}$ such that $\varphi\in
  C^\infty (\U_{-k},\R^{2n})$
\item[2)] $e^{-q^j(u)JA_j}u-\eta_{p(u)}\in\Pi_{p(u)} \V^{-k}$.
\end{itemize}
\end{lemma}
\proof Consider the condition 2). It is equivalent to the couple of
equations
\begin{align}
\label{e.con.coo}
0=f_l(q,p,u):=\langle e^{-q^jJA_j}u-\eta_p;A_l\eta_p\rangle\equiv
\langle u;e^{q^jJA_j}A_l\eta_p \rangle-2p^j=0\ ,
\\
0=g^l(q,p,u):= \langle e^{-q^jJA_j}u-\eta_p;E\dep{\eta_p}l
\rangle\equiv \langle u; e^{q^jJA_j}E\dep{\eta_p}l  \rangle
-\langle\eta_p ; E\dep{\eta_p}l\rangle 
\end{align} 
Both the functions $f$ and $g$ are smoothing, so we try to apply the
implicit function theorem in order to define the functions $q(u)$,
$p(u)$. First remark that the equations are fulfilled at
$(q,p,u)=(0,p_0,\eta_{p_0})$, then we compute the derivatives of such
functions with respect to $q^j,p_j$ and show that they are
invertible. We have
\begin{eqnarray*}
\dep{f_j}k\big\vert_{(0,p_0,\eta_{p_0})}=\left[\langle u;e^{q^l
    JA_l}A_j\dep{\eta
    _{p}}k\rangle-2\delta_j^k\right]_{(0,p_0,\eta_{p_0}))}=
-\delta_j^k\ , 
\end{eqnarray*}
where we used 
\begin{equation}
\label{e.con.1}
\delta_j^k=\dep{\null}k\frac{1}{2}\langle \eta_{p};A_j\eta_{p}\rangle=
\langle \eta_{p}; A_j \dep{\eta_p}k\rangle\ . 
\end{equation}
Then we have
\begin{equation}
\label{e.eq.q}
\frac{\partial f_j}{\partial
  q_k}\big\vert_{(0,p_0,\eta_{p_0})}=\langle\eta_{p_0};
JA_jA_k\eta_{p_0} \rangle=0 
\end{equation}
by the skew-symmetry of $J$ and property (S1).

We come to $g$. 
\begin{equation}
\label{e.eq.1}
\dep{g^j}k\big\vert_{(0,p_0,\eta_{p_0})}=\langle\eta_{p_0};
E\frac{\partial^2\eta_{p_0} 
}{\partial p_j\partial p_k}\rangle-\langle\dep {\eta_{p_0}}k;E\dep{\eta_{p_0}}j
\rangle   -\langle\eta_{p_0}; 
E\frac{\partial^2\eta_{p_0}  
}{\partial p_j\partial p_k}\rangle  
\end{equation} 
which vanishes by (H6). Finally we have 
$$
\frac{\partial g^j}{\partial q^k}\big\vert_{(0,p_0,\eta_{p_0})}=\langle
A_k\eta_{p_0};\dep {\eta_{p_0}}j \rangle=\delta_k^j 
$$
Therefore the implicit function theorem applies and gives the result.
\qed

\begin{corollary}
\label{cor.srt}
Any function $u\in \Hc^k$ in a neighborhood of $\eta_{p_0}$ can be
uniquely represented as 
\begin{equation}
\label{eq.rep.q}
u=e^{q^jJA_j}\left(\eta_{p}+\pip\phi  \right)\ ,
\end{equation}
with $(q,p,\phi)\in \R^n\times \R^n\times \V^k$ smoothly dependent on
$u$.
\end{corollary}

Remark that $p(u)\not= \Ph(u)$.

%
%
%

\begin{lemma}
\label{l.red.1}
Fix $\phi\in \V^l$, and let $X\in T_{i(\phi)}\Sc\cap \Hc^k$. Assume
$k+l\geq d_A$, then there exist $Q\equiv (Q^1,...,Q^n)\in\R^n$ and
$\Phi\in \Hc^{\min\{k,l-d_A\}}$ such that
\begin{equation}
\label{e.red.4}
X=Q^j JA\relax_j i(\phi)+i_*\Phi\ .
\end{equation}
\end{lemma}
\proof We write explicitly the formula \eqref{e.red.4} and show how
to solve it for $Q$ and $\Phi$:
\begin{equation}
\label{e.red.5}
X=Q^jJA\relax_j(\eta_p+\pip \phi)+\left(\dep{\eta_p}j+\depp
j\phi\right)\langle \nabla p\relax_j;\Phi\rangle+\pip \Phi\ .
\end{equation} 
 Apply to such a formula the operator $\widetilde
{\pip}^{-1}\pip$ (recall that $\widetilde \Pi_p:\Pi_p \Hc^i\to\V^i$ is an
isomorphism) getting
\begin{equation}
\label{e.red.6}
\Phi=\widetilde{\pip}^{-1}\pip X- Q^j\widetilde {\pip}^{-1}\pip
JA\relax_j(\eta_p+\pip \phi) -\widetilde{\pip}^{-1}\pip
\left(\dep{\eta_p}j+\depp j\phi\right)\langle \nabla
p_j;\Phi\rangle\ .
\end{equation} 
Taking the scalar product with $\nabla p\relax_k$ we get 

\begin{align}
\nonumber
\left\langle \nabla p\relax_k;\Phi \right\rangle=
\left\langle \nabla p\relax_k ; \widetilde{\pip}^{-1}\pip X\right\rangle -
Q^j \left\langle \nabla p\relax_k ;\widetilde {\pip}^{-1}\pip
JA\relax_j(\eta_p+\pip \phi)\right\rangle 
\\
\label{e.red.7}
-\left\langle \nabla p\relax_k
;\widetilde{\pip}^{-1}\pip \left(\dep{\eta_p}j+\depp
j\phi\right)\right\rangle \langle \nabla p\relax_j;\Phi\rangle
\end{align}
By formula \eqref{smo.01} and by remark \ref{r.5.1} one has 
$$
\left\langle \nabla p\relax_k ; \widetilde{\pip}\relax_{-1}\pip X\right\rangle=
- \sum_l M^l_k\langle A\relax_l\phi; X\rangle+{\rm smoothing\ function}\ , 
$$
which is well defined under the assumptions of the lemma. Now, 
\begin{eqnarray*}
\left\langle \nabla p_k ;\widetilde {\pip}^{-1}\pip
JA\relax_j\pip \phi\right\rangle=
- \sum_l M^l_k\langle A\relax_l\phi; JA\relax_j\phi \rangle+{\rm
  smoothing\ function}
\\
= {\rm smoothing\ function}\ .
\end{eqnarray*}

Furthermore the coefficient of $ \langle \nabla p\relax_j;\Phi\rangle$
at r.h.s.~of \eqref{e.red.7} is small (and smoothing) if $\phi$ is small enough. Thus one
can solve \eqref{e.red.7} and compute $ \langle \nabla
p\relax_j;\Phi\rangle$ as a function of well defined objects and of
$Q^j$.

Take the scalar product of \eqref{e.red.5} with
$E\dep{\eta_p}l,$ getting
\begin{eqnarray*}
\langle X;E\dep {\eta_p}l\rangle=Q^j\left(\langle E\dep {\eta_p}l
;JA\relax_j\eta_p\rangle + \langle E\dep {\eta_p}l ; JA\relax_j\pip\phi\right)
\\ + \left( \langle E\dep {\eta_p}l;\dep {\eta_p}j\rangle +\langle
E\dep {\eta_p}l;\depp j\phi\rangle \right)\langle\nabla
p\relax_j;\Phi\rangle
\\
= Q^j(-\delta^l\relax_j-\langle\dep{\eta_p}l;A\relax_j\pip\phi\rangle
)+\langle\nabla p\relax_j;\Phi\rangle \langle
E\dep {\eta_p}l;\depp j\phi\rangle\ .
\end{eqnarray*}
Substitute the expression we got for $\langle\nabla p\relax_j;\Phi\rangle$,
and then it is immediate to see that it is possible to compute the
$Q_j$'s, and thus also $\langle \nabla p_j;\Phi\rangle$, and use
\eqref{e.red.6} to get $\Phi$. \qed

\begin{lemma}
\label{l.reg.2}
Take $\phi\in\V^k$ with $k$ large enough, then the following formula
holds
\begin{equation}
\label{e.red.7.1}
X_H(e^{q^jJA\relax_j}i(\phi)
)=Q^l(\phi)e^{q^jJA\relax_j}JA\relax_li(\phi)+e^{q^jJA\relax_j}i_*X_{H_r}(\phi)\ .
\end{equation}
Furthermore there exists a matrix $\tilde M^l_j=\delta^l\relax_j+\hat
M^l\relax_j$ with $\hat M^l\relax_j$ smoothing functions, such that
\begin{equation}
\label{e.red.7a}
Q^l=\tilde M\relax_j^l\di H \dep{\eta_{p_0}}j\ .
\end{equation}
\end{lemma}
\proof First, the vector field $X_H$ is equivariant,
i.e. $$X_H(e^{q^jJA\relax_j}\phi)=e^{q^jJA\relax_j}X_{\phi}(\phi)\ ,$$
thus it is enough to verify the formula for $q^j=0$. Furthermore
$X_H(i(\phi))\in T_{i(\phi)}\Sc$ thus, by the preceding lemma it
admits the representation
\begin{equation}
\label{e.red.8}
X_H=Q^lJA\relax_li(\phi)+i_*\Phi\ ,
\end{equation}
and remark that, for any choice of $l$, one has
$$
\omega(JA\relax_li(\phi);i_*\Psi )= \langle A\relax_li(\phi);i_*\Psi
\rangle=0\ ,\quad \forall \Psi\in \V 
$$
since this is the condition ensuring that $i_*\Psi\in
T_{i(\phi)}\Sc$. 

Remark that (by lemma \ref{l.red.1}), at such points, any vector $U\in \Hc^k$
admits the representation
\begin{equation}
\label{e.red.10}
U=\alpha^lJA\relax_li(\phi)+i_*\Psi+\beta^j\dep {\eta_{p}}j\ ;
\end{equation}
we insert such a representations in the definition of the vector field
$X_H$. Obtaining
\begin{eqnarray*}
\di H U= \alpha^l \di H JA\relax_l i(\phi)+\di Hi_*\Psi +\beta^j \di H\dep
    {\eta_{p}}j 
\\
= \omega(X_H;U)= \omega(Q^lJA\relax_li(\phi)+i_*\Phi;
\alpha^lJA\relax_li(\phi)+i_*\Psi+\beta^j\dep {\eta_{p}}j ) \ ,
\end{eqnarray*}
which, exploiting the invariance of $H$ and \eqref{e.red.10}, gives 
\begin{equation}
\label{e.red.11}
\di Hi_*\Psi +\beta^j \di H\dep
    {\eta_{p}}j=Q^l\beta_j\omega(JA\relax_li(\phi);\dep {\eta_{p}}j)+
    \omega(i_*\Phi;i_*\Psi) +\beta^j \omega(i_*\Phi;\dep
          {\eta_{p}}j)\ .
\end{equation}
Taking $\beta^j=0$ we get $\di(i^*H)\Psi=i^{*}\omega(\Phi;\Psi)$, which
shows that $\Phi=X_{H_r}$. To get the formula for  the $Q$'s take
$\Phi=0$ and all the $\beta$'s equal to zero but one. Thus we get
\begin{eqnarray*}
\di H\dep{\eta_p}j=Q^l\omega(JA\relax_li(\phi);\dep{\eta_p}j)= Q^l\langle
A\relax_l(\eta_p+\pip \phi);\dep {\eta_p}j\rangle
\\
= Q^l \left(\delta^j\relax_l+\langle
A\relax_l\pip \phi;\dep {\eta_p}j\rangle  \right)
\end{eqnarray*}
which gives the thesis.\qed

From this Lemma the thesis of the theorem immediately follows.

\section{Proof of Perelman's Lemma \ref{perelman}}\label{B}

First we transform the equation $\dot \phi_c=L(t)\phi_c$ to a more
suitable form. To this end we make the transformation
\begin{equation}
\label{tilde}
\phi=e^{q^j(t)JA_j}\tilde \phi\ ,\quad \dot q^j=w^j\ ,\quad q^j(0)=0\ ,
\end{equation}
which gives
\begin{equation}
\label{tilde.1}
\frac{d}{dt}\tilde \phi =P_c(t)H(t)\tilde \phi
-\tilde R\tilde \phi\ ,
\end{equation}
where 
\begin{equation}
\label{tilde.11}
H(t):= L_{00}+J\tilde V(t)
\end{equation}
and
\begin{align*}
L_{00}&:=J(A_0+\E A_4)
\\
P_c(t)&:=e^{-q^JA_j}P_ce^{q^jJA_j} \ ,\quad \tilde
V(t):=e^{-q^jJA_j}V_0e^{q^jJA_j}\ ,\quad 
\\
\tilde
R&:=w^j\left[P_c(t)-\uno\right]
JA_jP_c(t)+w^jP_c(t)JA_j\left[P_c(t)-\uno\right] \ . 
\end{align*}
Explicitly the operators $\tilde V(t)$ and $P_c(t)$
can be computed by remarking that, since $e^{q^jJA_j}$ is canonical
and unitary for any fixed time, one has
\begin{eqnarray*}
d^2H_P(b_\E)(e^{q^jJA_j}\tilde \phi,e^{q^jJA_j}\tilde \phi)=
\frac{1}{2}\left\langle V_0e^{q^jJA_j}\tilde \phi;e^{q^jJA_j}\tilde
\phi \right\rangle =\frac{1}2\langle \tilde V \tilde \phi;\tilde \phi
\rangle
\\
=d^2H_P(e^{-q^jJA_j}b_\E)(\tilde \phi,\tilde \phi)\ .
\end{eqnarray*}
Thus the projector $P_c(t)$ is the projector on the continuous
spectrum of $H(t)$. From this and the fact that
$\e^{q^jJA_j}\tilde\phi\in P_c\V$, it follows in particular that, for
any time $t$, one has
$P_c(t)\tilde \phi(t)=\tilde \phi(t)$.
Remark also that one has 
\begin{align}
\label{pd}
\tilde R= w^j\left[P_c(t)-\uno\right]
JA_jP_c(t)+w^jP_c(t)JA_j\left[P_c(t)-\uno\right]
\\
\nonumber
= 
-w^j[P_d(t)JA_jP_c(t)+P_c(t)JA_jP_d(t)] 
\end{align}
which therefore is a small smoothing operator.
{\bf Omitting tildes} one has the explicit formula
\begin{align}
\label{vdit}
JV(t)\phi=-\beta'(b_{\E}^2(x-\bq(t)))\phi
\\
\nonumber
-\beta''( b_{\E}^2(x-\bq(t))
)2\Re\left( \e ^{-\im q^4(t)}b_{\E}(x-\bq(t))\phi \right) \e^{\im
  q^4(t)}b_{\E}(x-\bq(t))\ .
\end{align}
Here and below we denote by
$\bq\in\R^3$ the vector with components $q^j$, $j=1,2,3$.

We work on the
equation 
\begin{equation}
\label{tilde.12}
\dot \phi=H(t)\phi+R(t)\phi\ ,
\end{equation} 
following almost literally the proof given by Perelman. With a slight
abuse of notation we will here denote by $\U(t,s)$ the evolution
operator of such an equation. Remark that, fro the
fact that $L(t)$ leaves $P_c\V^k$ invariant, one has
\begin{equation}
\label{pc.1}
P_c(t)\U(t,s)=\U(t,s)P_c(s)\ .
\end{equation}

First we have the following proposition (which follows from proposition 1.1
of \cite{Per04} and the remark that $\e^{q^jJA_j}$ conserves all the
$L^p_x$ norms)

\begin{proposition}
\label{per04}
There exists $\epsilon_0$ such that, provided $|w^j|=|\dot
q^j(t)|<\epsilon_0$, then one has
\begin{equation}
\label{per04.eq}
\sup_{a\in\R^3,\ t\in\R, \ s\in\R}\left(\norma{\langle
  x-a\rangle^{-\nu} \e^{H(t)s}P_c(t)\phi}_{L^2}\langle s\rangle^{3/2}
\right) \sleq (\norma{\phi}_{L^2}+\norma{\phi}_{L^1})\ .
\end{equation}
\end{proposition}
We are going to prove the following local decay estimate from
which the Strichartz type inequalities \eqref{flow} and \eqref{flow1}
follow.
\begin{lemma}
\label{loc.dec}
The evolution operator $\U$ satisfies
\begin{align}
\label{loc.eq.1}
\norma{\langle
  x-a\rangle^{-\nu} \U(t,s)P_c(s)\phi}_{L^2} \sleq
\frac{\norma{\phi}_{L^2}+\norma{\phi}_{L^1}}{\langle
  t-s\rangle^{3/2}}\ ,
\\
\nonumber
\forall t\geq s\ ,\quad \forall a\in\R^3\ .
\end{align}
\end{lemma}
\proof It is clearly sufficient to work with $s=0$. We will make use
of the following Duhamel formula 
\begin{eqnarray*}
\phi(t)=\e ^{H(t)t}\phi_0+ \int_0^t ds \e
^{H(t)(t-s)}[H(s)-H(t)]\phi(s)+\int_0^t ds \e
^{H(t)(t-s)}R(s)\phi(s)\ .
\end{eqnarray*}
Applying $P_c(t)$ and iterating once the formula one gets
$\phi(t)=I_1+I_2+I_3+I_4+I_5$,
where
\begin{eqnarray*}
I_1&=&\e ^{H(t)t}P_c(t)\phi_0
\quad I_2=\int_0^t ds\, \e
^{H(t)(t-s)}P_c(t)R(s)\phi(s)\ 
\\
I_3&=& \int_0^t ds\, \e
^{H(t)(t-s)}P_c(t)[H(s)-H(t)]\e^{H(s)s}P_c(s)\phi_0\ ,
\\
I_4&=& \int_0^t ds\int_0^sds_1 \e
^{H(t)(t-s)}P_c(t)[H(s)-H(t)]\e^{H(s)(s-s_1)}P_c(s)R(s_1)\phi(s_1)\ ,
\\
I_5&=& \int_0^t ds \int_0^sds_1\e
 ^{H(t)(t-s)}P_c(t)[H(s)-H(t)]\e^{H(s)(s-s_1)}P_c(s)[H(s_1)-H(s)]\phi(s_1)\ .
\end{eqnarray*}
The only nontrivial estimate is that of $I_5$. We start by the others
and then we concentrate on $I_5$,

The estimate of $I_1$ is an immediate consequence of proposition
\ref{per04}.  Define
\begin{equation}
\label{mm.1}
m(t):=\sup_{a\in\R^3,\ 0\leq \tau\leq t} \left(\norma{\langle
  x-a\rangle^{-\nu} \phi(t)}_{L^2}\langle\tau\rangle^{3/2}\right)\ ,
\end{equation} 
and, in order to estimate $I_2$ remark that, for $p=1,2$, one has
$$
\norma{R(s)\phi(s)}_{L^p}\leq \norma{\langle x-\bq(s)\rangle^{N}
  R(s)\phi(s)}_{L^2} \sleq \epsilon \norma{\langle
  x-\bq(s)\rangle^{-\nu} \phi(s)}_{L^2} \leq \epsilon
\frac{m(s)}{\langle s\rangle^{3/2}}
$$
Substituting in $I_2$ one gets
\begin{eqnarray*}
\norma{\langle x-a\rangle^{-\nu}I_2(t)}_{L^2}\sleq \int_0^t
ds\frac{1}{\langle
  t-s\rangle^{3/2}}\left(\norma{R(s)\phi(s)}_{L^1}+\norma{R(s)\phi(s)}_{L^2}
\right)
\\
\sleq \epsilon m(t)\int_0^t ds\frac{1}{\langle
  t-s\rangle^{3/2}}\frac{1}{\langle
  s\rangle^{3/2}}
\sleq m(t)\frac{1}{\langle
  t\rangle^{3/2}}
\end{eqnarray*}
The estimate of $I_4$ is similar. For estimating $I_3$, first remark
that $H(s)-H(t)=V(s)-V(t)$, and thus, by \eqref{vdit}, for any
function $\phi$, one has
$$ \left|\left[H(s)-H(t)\right]\phi\right|\sleq \left|\langle
x-\bq(s)\rangle^{-N} \phi \right|+\left|\langle
x-\bq(t)\rangle^{-N} \phi \right|\ ,\ \forall N\ .
$$
So, we have
\begin{eqnarray*}
&&\norma{\langle x-a\rangle^{-\nu}I_3(t)}_{L^2}\\
&&\sleq \int_0^t
\frac{ds}{\langle t-s\rangle^{3/2}} \left[ \norma{\langle
x-\bq(s)\rangle^{-\nu}\e^{H(s)s}P_c(s)\phi_0}_{L^2}+  \norma{\langle
x-\bq(t)\rangle^{-\nu}\e^{H(s)s}P_c(s)\phi_0}_{L^2}   \right]  
\\
&&\sleq  \int_0^t\frac{ds}{\langle t-s\rangle^{3/2} \langle
  s\rangle^{3/2}} (\norma{\phi_0}_{L^1}+\norma{\phi_0}_{L^2}) \sleq
\frac{1}{\langle t\rangle^{3/2}}
(\norma{\phi_0}_{L^1}+\norma{\phi_0}_{L^2})\ .
\end{eqnarray*}
We come to $I_5$. Here the key remark is that (with a slight abuse of
notation) 
\begin{equation}
\label{vtmeno}
|V(s_1)-V(s_2)|\sleq \epsilon^{1/2}\left\langle
x-\bq(s_1)\right\rangle^{-N} \ ,\quad \forall N 
\end{equation} 
and all $s_1$, $s_2$ with
$|s_1-s_2|\sleq \epsilon^{-1/2}$. 

Consider the two cases $t\leq 4\epsilon^{-1/2}$ and $t\geq
4\epsilon^{-1/2}$. Exploiting \eqref{vtmeno} one easily gets that in
the first case 
$$
\norma{\langle x-a\rangle^{-\nu}I_5(t)}_{L^2}\sleq
\epsilon^{1/2}\frac{m(t)}{\langle t\rangle^{3/2}}\ .
$$
In the second case  $t\geq
4\epsilon^{-1/2}$, split the interval of integration of $s$ into three
parts, accordingly define 
\begin{eqnarray*}
I_{51}=\int_0^{\epsilon^{-1/2}}ds\ ,\ I_{52}=\int_{\epsilon^{-1/2}}^{t-\epsilon^{-1/2}}
ds\ ,\ I_{53}=\int_{t-\epsilon^{-1/2}}^tds\ .\ 
\end{eqnarray*}
The term $I_{51}$ is estimated exploiting the fact that in the
considered interval $V(s)-V(s_1)$ fulfill the estimate
\eqref{vtmeno}. Thus one gets
$$
\norma{\langle x-a\rangle^{-\nu}I_{51}(t)}_{L^2}\sleq
\epsilon^{1/2}\frac{m(t)}{\langle t\rangle^{3/2}}\ .
$$ Similarly $I_{53}$ is estimated using the fact that in such an
interval $V(s)-V(t)$ fulfill the estimate \eqref{vtmeno}, and thus it
gives the same contribution as $I_{51}$. Finally concerning $I_{52}$,
one has
\begin{eqnarray*}
\norma{\langle x-a\rangle^{-\nu}I_{52}(t)}_{L^2}\sleq
=\int_{\epsilon^{-1/2}}^{t-\epsilon^{-1/2}} ds \int_0^s ds_1
\frac{1}{\langle t-s \rangle^{3/2}}\frac{1}{\langle s-s_1
  \rangle^{3/2}} \frac{m(t)}{\langle s_1 \rangle^{3/2}} \\ \sleq m(t)
\int_{\epsilon^{-1/2}}^{t-\epsilon^{-1/2}}ds \frac{1}{\langle t-s
  \rangle^{3/2}}\frac{1}{\langle s \rangle^{3/2}}\sleq
\frac{\epsilon^{1/4}m(t)}{\langle t\rangle^{3/2}}\ .
\end{eqnarray*}
Collecting all the results one gets
\begin{equation}
\label{Ifin}
m(t)=\sup \left(\langle t\rangle^{3/2} \norma{\langle
  x-a\rangle^{-\nu}\phi(t)}_{L^2} \right)\sleq \norma{\phi_0}_{L^2}+
\norma{\phi_0}_{L^1}+\epsilon^{1/4}m(t)\ , 
\end{equation}
from which the thesis immediately follows. \qed

\noindent {\it End of the proof of lemma \ref{perelman}.} Consider the
following Duhamel formulae
\begin{align}
\label{per.222}
\U(t,0)&P_c(0)\phi_0=\e^{tL_{00}}P_c(0)\phi_0
\\
\nonumber
&+\int_0^t \di
s\e^{L_{00}(t-s)}(V(s)+R(s))\U(s,0)P_c(0)\phi_0
\\
\U(t,0)&P_c(0)\phi_0=P_c(t)\e^{tL_{00}}P_c(0)\phi_0
\\
\nonumber
&+\int_0^t \di
s\U(t,s)P_c(s)(V(s)+R(s))\e^{L_{00}s}\phi_0\ .
\end{align}
Inserting the second one in the integral of the first one one gets
that the quantity to be estimated is the sum of three integrals,
which can be easily estimated using \eqref{loc.eq.1} and the fact that
$\e^{L_ct}$ fulfills the Strichartz estimate \eqref{flow} as proved
e.g. in \cite{Cuc01} or \cite{Per04}.

The retarded estimate \eqref{flow1} can be deduced from \eqref{flow}
by reproducing exactly the argument by Keel and Tao.

The verification of (St.3\null) is a small variant and is omitted.
\qed



\begin{thebibliography}{FGJS04}

\bibitem[BC11]{BC09}
D.~Bambusi and S.~Cuccagna, \emph{On dispersion of small energy solutions of
  the nonlinear {K}lein {G}ordon equation with a potential},
Amer. J. Math.  \textbf{133}
  (2011), no.~5, 1421-1468.

\bibitem[Bec12]{Bec11}
Marius Beceanu, \emph{A critical center-stable manifold for {S}chr\"odinger's
  equation in three dimensions}, Comm. Pure Appl. Math. \textbf{65} (2012),
  no.~4, 431--507. 

\bibitem[BP92]{BusPer92}
V.~S. Buslaev and G.~S. Perelman, \emph{Scattering for the nonlinear
  {S}chr\"odinger equation: states that are close to a soliton}, Algebra i
  Analiz \textbf{4} (1992), no.~6, 63--102.

\bibitem[CM08]{CucMiz}
S.~Cuccagna and T.~Mizumachi, \emph{On asymptotic stability in energy space of
  ground states for nonlinear {S}chr\"odinger equations}, Comm. Math. Phys.
  \textbf{284} (2008), no.~1, 51--77.

\bibitem[Cuc01]{Cuc01}
S.~Cuccagna, \emph{Stabilization of solutions to nonlinear {S}chr\"odinger
  equations}, Comm. Pure Appl. Math. \textbf{54} (2001), no.~9, 1110--1145.

\bibitem[Cuc11a]{Cuc09}
\bysame, \emph{The {H}amiltonian structure of the nonlinear {S}chr\"odinger
  equation and the asymptotic stability of its ground states}, Comm. Math.
  Phys. \textbf{305} (2011), no.~2, 279--331.

\bibitem[Cuc11b]{CucMov}
\bysame, \emph{{O}n asymptotic stability of moving ground states of the
  nonlinear {S}chr\"odinger equation}, Preprint: http://arxiv.org/abs/1107.4954
  (2011).

\bibitem[FGJS04]{FroJon}
J.~Fr{\"o}hlich, S.~Gustafson, B.~L.~G. Jonsson, and I.~M. Sigal,
  \emph{Solitary wave dynamics in an external potential}, Comm. Math. Phys.
  \textbf{250} (2004), no.~3, 613--642.

\bibitem[GNT04]{GSNT04}
Stephen Gustafson, Kenji Nakanishi, and Tai-Peng Tsai, \emph{Asymptotic
  stability and completeness in the energy space for nonlinear {S}chr\"odinger
  equations with small solitary waves}, Int. Math. Res. Not. (2004), no.~66,
  3559--3584.

\bibitem[GS07]{GS07}
Zhou Gang and I.~M. Sigal, \emph{Relaxation of solitons in nonlinear
  {S}chr\"odinger equations with potential}, Adv. Math. \textbf{216} (2007),
  no.~2, 443--490.

\bibitem[GW08]{GW08}
Zhou Gang and M.~I. Weinstein, \emph{Dynamics of nonlinear
  {S}chr\"odinger/{G}ross-{P}itaevskii equations: mass transfer in systems with
  solitons and degenerate neutral modes}, Anal. PDE \textbf{1} (2008), no.~3,
  267--322.

\bibitem[Per04]{Per04}
G.~S. Perelman, \emph{Asymptotic stability of multi-soliton solutions for
  nonlinear {S}chr\"odinger equations}, Comm. Partial Differential Equations
  \textbf{29} (2004), no.~7-8, 1051--1095.

\bibitem[Per11]{appPer}
\bysame, \emph{Asymptotic stability in {$H^1$} of {NLS}. {O}ne soliton case},
  Personal communication (2011).

\bibitem[Sch87]{Schmid}
R.~Schmid, \emph{Infinite-dimensional {H}amiltonian systems}, Monographs and
  Textbooks in Physical Science. Lecture Notes, vol.~3, Bibliopolis, Naples,
  1987.

\bibitem[Sig93]{Sig93}
I.~M. Sigal, \emph{Nonlinear wave and {S}chr\"odinger equations. {I}.
  {I}nstability of periodic and quasiperiodic solutions}, Comm. Math. Phys.
  \textbf{153} (1993), no.~2, 297--320.

\bibitem[SW99]{SofW99}
A.~Soffer and M.~I. Weinstein, \emph{Resonances, radiation damping and
  instability in {H}amiltonian nonlinear wave equations}, Invent. Math.
  \textbf{136} (1999), no.~1, 9--74.

\bibitem[Yaj95]{Yaj}
K.Yajima, {\em The $W^{k,p}$ continuity of wave operators for
Schr\"odinger operators \/},  J. Math. Soc. Japan, 47
 (1995), pp.
  551--581.

\end{thebibliography}

\providecommand{\bysame}{\leavevmode\hbox to3em{\hrulefill}\thinspace}
\providecommand{\MR}{\relax\ifhmode\unskip\space\fi MR }
\providecommand{\MRhref}[2]{%
  \href{http://www.ams.org/mathscinet-getitem?mr=#1}{#2}
}
\providecommand{\href}[2]{#2}

\noindent 
D.~Bambusi: Dipartimento di Matematica ``Federico Enriques'',
Universit\`a degli Studi di Milano, Via Saldini 50, 20133 Milano,
Italy.

\noindent  {\it E-mail}: {\tt dario.bambusi@unimi.it}

\end{document}